\documentstyle[12pt,leqno]{article}
\title{Foliations with Transversal Quaternionic Structures}
\author{{\normalsize by}\\Paolo Piccinni and Izu Vaisman}
\date{}
\begin{document}
\maketitle
{\def\thefootnote{*}\footnotetext[1]%
{{\it 1991 Mathematics Subject Classification}
53 C 12, 53 C 25. \newline\indent{\it Key words and phrases}:
Foliation,Quaternionic Structure, Twistor Space.\newline\indent
{\it Acknowledgements.} This paper was started while the second author
visited the Istituto di Matematica, Universit\`a di Roma 1, under the
auspices of the CNR, Italy. He is expressing here his thanks to the hosts.
The first author acknowledges support by MURST of Italy and support and
hospitality by the Erwin Schr\"odinger Institute for Mathematical
Physics in Vienna.}
\begin{center} \begin{minipage}{12cm}
A{\footnotesize BSTRACT.
We consider manifolds equipped with a foliation
$\cal F$ of codimension $4q$,
and an almost quaternionic structure $Q$ on the
transversal bundle of ${\cal F}$. After discussing conditions
of projectability and integrability of $Q$, we study the
transversal twistor space
$Z{\cal F}$ which, by definition, consists
of the $Q$-compatible almost complex
structures. We show that
$Z{\cal F}$ can be endowed with a lifted foliation
${\widehat {\cal F}}$ and two natural
almost complex structures $J_1$, $J_2$ on the transversal bundle of
$\widehat{\cal F}$.
We establish the conditions which ensure the
projectability of $J_1$ and $J_2$, and the integrability
of $J_{1}$ ($J_{2}$ is never integrable).}
\end{minipage}
\end{center}
\section{Preliminaries}
We recall the basic definitions of quaternionic geometry e.g.,
\cite{{B},{Sal}}. The general framework is the $C^{\infty}$ category.

An {\em almost hypercomplex structure} on a $4q$-dimensional differentiable
manifold $N^{4q}$ is an ordered
triple $H=(I_{1},I_{2},I_{3})$ of almost complex
structures satisfying the quaternionic identities $I_{\alpha}\circ
I_{\beta}=I_{\gamma}$ for $(\alpha,\beta,\gamma)=(1,2,3)$ and cyclic
permutations. If the structures $I_{1},I_{2},I_{3}$ are integrable, $H$ is
said to be a {\em hypercomplex structure}.

If $H=(I_{1},I_{2},I_{3})$ is an almost hypercomplex structure, any triple
$(J_{1},J_{2},$
$J_{3})$ obtained from $(I_{1},I_{2},I_{3})$ by multiplying
by a matrix of $SO(3)$ is again an almost hypercomplex structure.
Moreover, there exists a set of {\em compatible almost complex structures}
associated with a given almost hypercomplex structure namely,
the set of all $J=a_{1}I_{1}+a_{2}I_{2}+a_{3}I_{3}$ where
$a_{1},a_{2},a_{3}$ are functions satisfying
$a_{1}^{2}+a_{2}^{2}+a_{3}^{2}=1$.

An {\em almost quaternionic structure} on the manifold $N^{4q}$ is a rank
$3$ vector subbundle $Q$ of the endomorphism bundle $End(TN)$ locally
spanned by almost hypercomplex structures $H=(I_{1},I_{2},I_{3})$ which are
related by $SO(3)$-matrices on the intersections of trivializing open sets.
A {\em quaternionic structure} on the manifold $N^{4q}$ is an
almost quaternionic structure such that there exists a torsionless
connection $\nabla$ of $TN$ which, when extended to the vector bundle
$End(TN)$, preserves the subbundle $Q$ i.e., $\nabla Q\subseteq Q$.
The existence of the $Q$-preserving torsionless connection $\nabla$
is not equivalent with the integrability of $Q$ as a $G$-structure.
The existence of a {\em flat}
torsionless connection which preserves $Q$ implies that $Q$ can be obtained
from local quaternionic coordinates on $N$.
If an almost quaternionic structure $Q$ is
fixed on $N^{4q}$, the {\em local bases}
$(I_{1},I_{2},I_{3})$ which span the vector bundle $Q$ are also called {\em
local compatible almost hypercomplex structures}, and any local
$J=a_{1}I_{1}+a_{2}I_{2}+a_{3}I_{3}$ with
$a_{1}^{2}+a_{2}^{2}+a_{3}^{2}=1$ is called a {\em local $Q$-compatible
almost complex structure}.

A Riemannian metric $g$ on a (almost) hypercomplex manifold $(N,H)$ is
{\em (almost)
hyperhermitian},
respectively {\em (almost) hyperk\"ahler}, if it is (almost)
Hermitian, respectively, (almost) K\"ahlerian,
with respect to all the structures
$I_{\alpha}$, $\alpha=1,2,3$, of $H$. (Then,
$g$ also is compatible with any $H$-compatible structure $J$.) Similarly, on
an almost quaternionic manifold $(N,Q)$, the metric $g$ is
{\em quaternion
Hermitian} if it is Hermitian with respect to the local bases $(I_{\alpha})$
of $Q$,
and it is {\em quaternion K\"ahler} if
it is quaternion Hermitian and $Q$ is parallel
(i.e., $\nabla Q\subseteq Q$) with respect to the
Levi-Civita connection $\nabla$ of $g$.
(In both cases, the property says nothing about the integrability
of the structures $I_{\alpha}$.)
The terms are also used for
manifolds endowed with the respective structures. Of course, a
hyperk\"ahler manifold necessarily is hypercomplex, and a quaternion
K\"ahler manifold necessarily is quaternionic.

The twistor space $ZN$ of an almost quaternionic manifold $(N,Q)$ is
defined as the manifold of the $Q$-compatible almost complex structures of the
tangent spaces of $N$.Thus, $ZN$ is an $S^{2}$-bundle associated
with the vector bundle $Q$, where $Q$ has the metric
which makes the local bases
$H=(I_{1},I_{2},I_{3})$ orthonormal bases \cite{{B},{Sal}}.

Now, let us consider a $C^{\infty}$ manifold $M^{p+4q}$, equipped with a
$p$-dimensional foliation ${\cal F}$. Denote by $L=T{\cal F}$ the tangent
bundle of ${\cal F}$, and by $\nu{\cal F}=TM/L$ its transversal vector
bundle of rank $4q$. We will often identify the transversal bundle
$\nu{\cal F}$ with a complementary distribution $E$ of $L$ i.e., a
splitting of the exact sequence
$$0\rightarrow L=T{\cal F}\stackrel
{\subseteq}{\rightarrow} TM\stackrel{\pi_{\nu}}{\rightarrow}
\nu{\cal F}\rightarrow0.$$

Almost hypercomplex and almost quaternionic structures
can be defined similarly on vector
bundles of rank $4q$. Accordingly, they will be reductions of the structure
group of the bundle to $G$, where the group $G=GL(q,{\bf H})$ for the
almost hypercomplex structures and
$$G=GL(q,{\bf H})\cdot Sp(1)=\frac{GL(q,{\bf H})\times Sp(1)}{\pm Id}
$$
for the almost quaternionic structures (${\bf H}$ is the algebra of the
quaternions). Furthermore, almost hyperhermitian and quaternion
Hermitian structures correspond to the structure groups
$Sp(q)$ and $Sp(q)\cdot Sp(1)$, respectively.

We will consider structures
of these types on the transversal bundle $\nu{\cal F}$ of a
foliation ${\cal F}$, and refer to them as {\em transversal almost
hypercomplex, transversal almost quaternionic, etc.
structures} of the foliation
${\cal F}$. (Such structures sporadically appeared in the literature e.g.,
\cite{NTV}.)

In what follows, we use {\em Bott connections} \cite{Bt}
$D:\Gamma TM\times\Gamma\nu{\cal F}\rightarrow\Gamma\nu{\cal F}.$
A Bott connection is a connection on the transversal bundle $\nu{\cal F}$
which extends the {\em partial connection}
$\stackrel{\circ}{D_{}}:\Gamma L\times\Gamma\nu{\cal F}
\rightarrow\Gamma\nu{\cal F}$
given by
$$\stackrel{\circ}{D_{}}_{Y}s=\pi_{\nu}[Y,X_{s}],\leqno{(1.1)}$$
where $Y$ is a tangent vector field of the leaves of the foliation ${\cal
F}$, and $X_{s}$ is any vector field on $M$
such that $\pi_{\nu}X_{s}=s$, $s\in\Gamma\nu{\cal F}$.
($\Gamma$ always denotes
spaces of global cross sections of vector bundles.)
Notice that an identification $\nu{\cal F}\approx E$, where
$TM=E\oplus L$, implies the replacement of (1.1) by
$$\stackrel{\circ}{D_{}}_{Y}X=\pi[Y,X],\hspace{5mm}X\in\Gamma
E,\leqno{(1.1')}$$ $\pi$ being the projection $\pi:TM\rightarrow E$.

A Riemannian metric $g$ splits $TM = T\cal F \oplus T^\perp \cal
F$, and we will take $E=T^\perp \cal
F \approx \nu \cal F$. Then, in particular,
$\stackrel{\circ}{D_{}}$ can be extended
to a Bott connection
$D$, by defining
$D_{X}=\pi\circ\nabla_{X}$ $(X\in\Gamma E)$,
where $\nabla$ is
the Levi-Civita connection of $g$.
For Riemannian foliations, this Bott
connection $D$ is the unique
torsionless metric connection of the normal bundle $T^\perp
\cal F \approx \nu \cal F$ \cite{Mol}.

\proclaim 1.1 Definition. {\rm (a)  A transversal almost hypercomplex structure
$H=(I_1,I_2,I_3)$ of $\cal F$
is} projectable {\rm if the partial connection
$\stackrel{\circ}{D_{}}$
preserves the structures
$I_\alpha$, i. e. $\stackrel{\circ}{D_{}} I_\alpha =0$, $\alpha=1,2,3$.}\\
{\rm (b)  A transversal almost quaternionic structure
$Q \subseteq End (\nu
\cal F)$} is projectable {\rm if $\stackrel{\circ}{D_{}}$ preserves $Q$:
$\stackrel{\circ}{D_{}} Q \subseteq Q$.}
\par
The projectability condition of $Q$ can be formulated
in terms of local bases
$H=(I_1,I_2,I_3)$. Namely, $Q$ is
projectable iff, for any choice of a local
basis
$H=(I_1,I_2,I_3)$,
there exist local 1-forms $\alpha, \beta, \gamma$ such that:
$$\stackrel{\circ}{D_{}}I_{1}=\alpha I_{2}+\beta I_{3},\;\;\;
\stackrel{\circ}{D_{}}I_{2}=-\alpha I_{1}+\gamma I_{3},\;\;\;
\stackrel{\circ}{D_{}}I_{3}=-\beta I_{1}-\gamma I_{2}.\leqno{(1.2)} $$
As a matter of fact, if equations of the type
(1.2) hold for some choice of $H$ similar equations hold for any choice of
$H$.

If $J\in\Gamma End(\nu{\cal F})$ we may also see it as a
cross section of $End\,E$, and for $Y \in \Gamma
T\cal F$, $s \in \Gamma \nu \cal F$ we have
$$(\stackrel{\circ}{D_Y} J)s = ~\stackrel{\circ}{D_Y} Js - J
\stackrel{\circ}{D_Y} s =
\pi_{\nu}[Y,JX_s]- J\pi_{\nu} [Y,X_s] \leqno{(1.3)}$$
for any $X_s \in \Gamma TM$ such that $s=\pi_{\nu} X_s$.
A cross section $s$ is
{\em projectable} if $s$ projects to a tangent vector field
of any local {\em space of slices} of the foliation ${\cal F}$ \cite{Mol}.
From (1.3), it follows that $J$ is projectable
iff $Js$ is projectable whenever $s$ is projectable.
Therefore, projectability in the sense of Definition 1.1 (a) means that we have
a structure which is the lift of almost hypercomplex structures of the
local slice spaces. The same is true in the case of Definition 1.1 (b) (see
Proposition 3.1 later on).

Accordingly, we will give
\proclaim 1.2 Definition. {\rm A
projectable, transversal, almost hypercomplex or almost quaternionic
structure of a foliation ${\cal F}$ is} integrable
{\rm if the projected structures
of the local slice spaces are hypercomplex or quaternionic, respectively.}
\par
If integrability holds, the word {\em almost} will be omitted.
\section{Examples}
We begin by an example of a foliation with a transversal almost
hypercomplex structure. Let ${\cal F}$ be a transversally holomorphic
foliation of real codimension $2q$ on a manifold $N^{p+2q}$. This means
that, on $N$, there are local coordinates $(y^{u},z^{a},\bar z^{a})$, where
$(y^{u})$ are real coordinates, and $(z^{a})$ are complex coordinates with
holomorphic transition functions, such that ${\cal F}$ is defined by
$z^{a}=const.,\,\bar z^{a}=const.$ Furthermore, assume that $g$ is a
bundle-like Riemannian metric on $N$, which is ${\cal F}$-transversally
K\"ahler i.e.,
$$g=g_{a\bar b}(z,\bar z)dz^{a}\otimes d\bar z^{b}
+g_{{\bar a}b}(z,\bar z)d\bar z^{a}\otimes dz^{b}+\cdots, \leqno{(2.1)}$$
where the first two terms define a K\"ahler metric in the coordinates
$(z)$, and the remaining (unexplicited) terms contain $dy^{u}$. (In this
paper, we use the Einstein summation convention.)

Now, consider the manifold $M$ defined by the total space of the conormal
bundle of the foliation ${\cal F}$ i.e., the annihilator of the tangent
bundle $T{\cal F}$.
If $E$ is the $g$-orthogonal
bundle of ${\cal F}$, then $M=E^{*}$. On $M$, there exists a natural lift
${\cal F}^{*}$ of ${\cal F}$ such that the leaves of ${\cal F}^{*}$ are
covering spaces of the leaves of ${\cal F}$. Moreover, the local slice
spaces of ${\cal F}^{*}$ are cotangent bundles of K\"ahler manifolds. It is
not difficult to construct an almost hypercomplex structure on such a
cotangent bundle. If we do this on the local slice spaces, the obtained
structures glue up to a projectable almost hypercomplex structure
transversal to ${\cal F}^{*}$. The construction of the almost hypercomplex
structure of the cotangent bundle of a K\"ahler manifold was described in
\cite{V7}. For the reader's convenience, we present here
the previously mentioned transversal structure of ${\cal F}^{*}$
directly.

As in \cite{{V1},{V4}}, let
$$\theta^{u}=dy^{u}+t_{a}^{u}dz^{a}+\bar t^{u}_{a}d\bar
z^{a}\leqno{(2.2)}$$
be a basis of the annihilator of $E$. Then, $\forall\zeta\in E^{*}$ we have
$$\zeta=\zeta_{a}dz^{a}+\bar\zeta_{a}d\bar z^{a},\leqno{(2.3)}$$
$(y^{u},z^{a},\bar z^{a},\zeta_{a},\bar\zeta_{a})$ are local coordinates
on $M$, and the system of equations
$$z^{a}=const.,\;\bar z^{a}=const.,\;
\zeta_{a}=const.,\;\bar\zeta_{a}=const.$$
defines the foliation ${\cal F}^{*}$. Obviously, ${\cal F}^{*}$ is again a
transversally holomorphic foliation, and we will denote by $I_{1}$ the
corresponding transversal complex structure. That is, $I_{1}$ is a complex
structure on the transversal bundle $\nu{\cal F}^{*}$ which, in turn,  may be
identified with the complementary bundle $S$ of $T{\cal F}^{*}$ given by
the equations $\theta^{u}=0$ on $M=E^{*}$.
As usual, we may identify $(S,I_{1})$ with the holomorphic part $S_{1,0}$
of $S\otimes_{\bf R}{\bf C}$.

Then, $S$ also has a canonical symplectic
structure namely, if seen on $M$, (2.3) is a $1$-form which may be viewed
as the ${\cal F}^{*}$-{\em transversal Liouville form} $\zeta$, and
$-d\zeta$ is the mentioned symplectic structure. When transferred to
$S_{1,0}$, these structures go to $\lambda=\zeta_{a}dz^{a}$ and $\omega=
dz^{a}\wedge d\zeta_{a}$, respectively.

Furthermore, the Levi-Civita connection of the transversal K\"ahlerian part
of $g$ yields a connection on $E$ with a corresponding horizontal
distribution ${\cal H}$ on $E^{*}$ given by \cite{V7}
$$d\zeta_{a}-\Gamma^{c}_{ab}\zeta_{c}dz^{b}=0,\;\;
d\bar\zeta_{a}-\bar\Gamma^{c}_{ab}\bar\zeta_{c}d\bar z^{b}=0,\leqno{(2.4)}$$
where the coefficients $\Gamma$ are the Christoffel symbols.
The equations (2.4) completed by $\theta^{u}=0$ define the {\em horizontal
part}
${\cal H}_{S}$ of the bundle $S$. Of course, $S$ is tangent to the fibers
of $E^{*}$ i.e., it contains the {\em vertical distribution}, say ${\cal
V}$, of this bundle. As a matter of fact, we have $S={\cal
H}_{S}\oplus{\cal V}$.

Now, continuing with the identification $(S,I_{1})\approx S_{1,0}$, we see
that a new complex structure $I_{2}$ of $S$ can be obtained by asking
$$I_{2}/_{{\cal H}_{S}}=\overline{\sharp_{\omega}\circ\flat_{g}},\;\;
I_{2}^{2}=-Id,\leqno{(2.5)}$$
where the {\em musical isomorphisms} are defined as in Riemannian
geometry and the bar denotes complex conjugation
(of course, only the transversal part of $g$ is used).

Finally, the same computations as in \cite{V7} show that
$(I_{1},I_{2},I_{3}:=I_{1}\circ I_{2})$ is an almost hypercomplex structure
on the vector bundle $S$, and this is the announced example

Now, we will describe two classes of examples of
foliations with projectable, transversal
quaternionic structure, which come from
3-Sasakian and quaternion Hermi\-tian-Weyl geometry, respectively
(cf. \cite{{BG},{BGM},{OP1},{OP2}}).
\par
A triple ($\xi^1$, $\xi^2$, $\xi^3$)
of orthonormal Killing vector fields on a $(4q+3)$-dimensional
Riemannian manifold $({\cal S},g)$ is said to define a 3-{\em
Sasakian structure}
if their brackets satisfy the identities $$[\xi^{\alpha},\xi^{\beta}]= 2
\xi^{\gamma}$$ ($(\alpha,\beta,\gamma)=(1,2,3)$ and cyclic permutations),
and, furthermore, the
dual 1-forms $\eta^\alpha=\flat_{g}\xi^{\alpha}$
satisfy the equations
$$(\nabla_Y \Phi^\alpha)Z=
\eta^\alpha(Z)Y-g(Y,Z)\xi^\alpha\leqno{(2.6)}$$
($\alpha=1,2,3$) where $\nabla$ is the
Levi-Civita connection of $g$, and
$\Phi^\alpha=\nabla \xi^\alpha\in End(T{\cal S})$.

A manifold $({\cal S},g)$ with a $3$-Sasakian structure is a $3$-{\em Sasakian
manifold}, and the vector fields $\xi^{1},\xi^{2},\xi^{3}$ span a
foliation ${\cal V}$ of ${\cal S}$. Furthermore, ${\cal V}$ is invariant by
the endomorphisms $\Phi^{\alpha}$, and it has the orthogonal distribution
$E$ defined by $\eta^1=0,\,\eta^2=0,\,\eta^3=0$. Therefore, $E$ also is
$\Phi^{\alpha}$-invariant. Following Proposition 1.2.4 of \cite{BG},
one has
$$({\Phi^\alpha})^2 = -I+ \eta^\alpha \otimes
\xi^\alpha,$$
and  $(I_1=- \Phi^1/_{E},\; I_2=- \Phi^2/_{E},\;
I_3=- \Phi^3/_{E})$ is an almost
hypercomplex structure on the distribution $E = T^\perp \cal F$.

Now, we will check that, although not every
$I_1, I_2, I_3$ is projectable, the vector bundle
$Q$ spanned by these structures is projectable
(see also \cite{{BG},{BGM}}) hence,
the foliation $\cal V$ has a projectable, transversal
quaternionic structure.

Let $X$ be a projectable cross section of $E$ i.e.,
$[\xi^{\alpha},X]\in\Gamma T{\cal V}$ for $\alpha=1,2,3$. Then
$$(\stackrel{\circ}{D_{}}_{\xi^{\alpha}}\Phi^{\beta})X=
\pi[\xi^{\alpha},\Phi^{\beta}X]=
\pi(\nabla_{\xi^{\alpha}}(\Phi^{\beta}X)-\nabla_{\Phi^{\beta}X}\xi^{\alpha})$$
$$=\pi((\nabla_{\xi^{\alpha}}\Phi^{\beta})X+\Phi^{\beta}(\nabla_{\xi^{\alpha}}X)
-\Phi^{\alpha}\circ\Phi^{\beta}X)$$
$$\stackrel{(2.6)}{=}\pi\{\Phi^{\beta}(\nabla_{X}\xi^{\alpha}+
[\xi^{\alpha},X])-\Phi^{\alpha}\circ\Phi^{\beta}X\}$$
$$=(\Phi^{\beta}\circ\Phi^{\alpha}-\Phi^{\alpha}\circ\Phi^{\beta})X
= 2(1-\delta_{\alpha\beta}) \Phi^\gamma X,$$
where if $\alpha\neq\beta$ then
$(\alpha, \beta, \gamma)$ is a cyclic permutations of
$(1,2,3)$.
The last equality holds because the structures $I^{\alpha}$ satisfy
the quaternionic identities.

We recall that compact 3-Sasakian manifolds ${\cal S}^{4q+3}$
where the foliation
$\cal V$ has all the leaves compact
project onto a compact positive quaternion K\"ahler
orbifold $N^{4q}$, and the leaves of $\cal V$
are homogeneous 3-dimensional spherical
space forms.
In the case of a regular foliation $\cal V$, the leaf space
$N^{4q}$ is
a positive quaternion K\"ahler manifold. Thus, the simplest example of a
foliation with projectable transversal quaternionic structure is the
Hopf fibration
$S^{4q+3} \rightarrow {\bf H}P^q$.
An example of a 3-Sasakian manifold where
$\cal V$ is not regular, but still
all the leaves are compact, is the following. Consider the
action of ${\bf Z}_3$ on the sphere $S^7 =\{(h_0,h_1)
\in {\bf H}^2\;/\;  h_0
\bar h_0 +h_1 \bar h_1 =1~\}$
generated by $(h_0,h_1) \mapsto (e^{\frac{2\pi i}{3}}h_0,
e^{\frac{4\pi i}{3}}h_1)$.
This action preserves the 3-Sasakian structure of $S^7$, therefore,
the quotient ${\bf Z}_3\backslash S^7$ is a 3-Sasakian manifold,
and its  foliation
$\cal V$ admits a projectable transversal quaternion K\"ahler structure.
In fact
this structure projects to the orbifold ${\bf Z}_3
\backslash{\bf H}P^1$ defined by
the induced action of $Z_{3}$. We refer the reader to  \cite{{BG},{BGM}} for
all these facts.

As a matter of fact, the projection on a quaternion K\"ahler manifold
always holds locally (cf. \cite{BG}, Theorem 2.3.4). This shows that the
transversal almost quaternionic structure of the foliation
${\cal V}$ of an arbitrary 3-Sasakian manifold always is
an integrable i.e., a quaternionic, structure.

A second class of examples of
foliations with a projectable transversal quaternionic structure is
that of the {\em locally conformal
quaternion K\"ahler manifolds} $M^{4q+4}$. This means that $M$ is endowed
with an almost
quaternionic structure $Q$ and a metric $g$,
which is Hermitian with respect to
the local compatible almost complex structures of $Q$,
and such that, over some open neighborhoods $\{U_i\}$ which cover $M$
($M=\cup_{i}U_{i}$), $g$ is
conformally related to local quaternion K\"ahler
metrics: $$g_{\vert_{U_i}} =e^{f_i}g'_i,$$
where $g'_i$ is
quaternion K\"ahler on $U_i$ and $f_{i}\in C^{\infty}(U_{i})$.

Such a structure defines the so called Lee 1-form
$\omega$, where $\omega_{\vert_{U_i}} = df_i$.
$\omega$ appears as a factor in the exterior
differential $d\Theta =
2\omega \wedge \Theta$ of the {\em K\"ahler $4$-form} $\Theta =
\sum_{\alpha=1}^{3}\Omega_{\alpha}\wedge
\Omega_{\alpha}$, where $\Omega_\alpha$ are the K\"ahler forms of the local
bases $(I_\alpha)$ ($\alpha=1,2,3$) of $Q$.
In the compact case and if
$g$ is not globally conformal quaternion K\"ahler,
a result of P. Gauduchon yields a metric in the conformal class
of $g$
such that its Lee form $\omega$ is parallel with respect
to the Levi-Civita connection
of the new metric \cite{{GAU},{OP1}}. With
this choice and the normalization
$\vert \omega \vert =1$,
the Lee vector field $\xi:=\sharp_{g}\omega$, and the
local vector fields $\xi^\alpha= I_\alpha \xi$
define a 4-dimensional foliation $\cal V$ (cf. \cite{OP1}, Proposition 1.7)
whose orthogonal bundle $E$ has a quaternionic structure $Q_{E}$
induced by the structure $Q$ of $M$ (again, see \cite{OP1}).

Moreover, the Lie derivative formulas of Proposition 1.7 of \cite{OP1}
allow for an easy verification of the fact that
$$\stackrel{\circ}{D_{}}_{\xi}(Q_{E})\subseteq Q_{E},
\stackrel{\circ}{D_{}}_{\xi^{\alpha}}(Q_{E})\subseteq Q_{E},$$
for any Bott connection $D$ of $E$.

Therefore, ${\cal V}$ is a foliation with a projectable transversal almost
quaternionic structure. Moreover, it follows from the proof of Theorem 5.1
of \cite{OP1} that this transversal structure is, in fact, integrable.

It is easy to
give examples where the foliation $\cal V$ is not a fibration over an
orbifold.
The simplest examples of locally conformal hyperk\"ahler
(hence, implicitly, quaternion K\"ahler) manifolds
are quotients $({\bf H}^2 - \{0\})/{\bf Z}$, where ${\bf Z}$
is an infinite cyclic group which
preserves  the
metric $g=(
h_0 \bar h_0 +h_1 \bar h_1)^{-1} g_0$,
conformal to the standard flat metric $g_0$.
If we take ${\bf Z}$
generated by $(h_0,h_1) \mapsto (2 h_0, 2 e^{{\sqrt
2}\pi i}h_1)$, we get
a foliation $\cal V$ with non compact leaves. Thus, no orbifold
structure is obtained on the
leaf space. However, the transversal quaternionic structure
of the foliation ${\cal V}$ defined above is still
projectable. (See \cite{OP1} for more explanations.)
\section{Projectability}
In this section we continue to use the notation of Section 1, and we discuss
the notion of projectability of a transversal almost
quaternionic structure
$Q \subset End (\nu \cal F)$ introduced by Definition 1.1.
\proclaim 3.1 Proposition. The almost quaternionic structure $Q \subset
End (\nu \cal F)$
is projectable iff $Q$ has
local compatible, projectable,
almost hypercomplex structures $(J_{\alpha})$ ($\alpha = 1,2,3$).
\par
\noindent{\bf Proof.}
Since local systems $(I_1, I_2, I_3)$, $(J_1, J_2, J_3)$ of $Q$-compatible
almost hypercomplex structures are $SO(3)$-related,
it follows that if $\stackrel{\circ}{D_{}}
J_\alpha =0$ then $\stackrel{\circ}{D_{}} I_\alpha$
are given by expressions of the type
(1.2), i. e. $Q$ is
projectable. Conversely,
since $\stackrel{\circ}{D_{}}$ is a flat partial connection
\cite{{Bt},{Mol}}, the
condition $\stackrel{\circ}{D_{}} Q \subset Q$
insures that $\stackrel{\circ}{D_{}}$ induces a flat
partial connection on the vector bundle $Q$. Accordingly,  frames
$J_1, J_2, J_3$ which are parallel with respect to $\stackrel{\circ}{D_{}}$
(i.e., $\stackrel{\circ}{D_{}}J_{\alpha}=0$) can be constructed.
Namely, if $V$ is a local transversal submanifold of ${\cal F}$,
we fix
$J_\alpha$ along $V$ then, translate them parallely along the local slices
of ${\cal F}$, with respect to an arbitrary Bott connection. Q.e.d.
\par If a decomposition $TM =E \oplus L$ is chosen,
$\nu \cal F$ is isomorphic with the subundle $E$ of $TM$, and the
following tensorial projectability criterion of an almost
complex structure $J$ of $E$
(i. e., $\stackrel{\circ}{D_{}} J = 0$) holds.
Let $\tilde J$ be the endomorphism of
$TM$ defined by
$$\tilde J (X) =\left\{ \begin{array}{ll}
JX, & \;\;\;\;{\rm for}\;X \in \Gamma E,  \\
0, & \;\;\;\;{\rm for}\;X \in \Gamma L,
\end{array}\right. \leqno{(3.1)}$$
and consider the {\em Nijenhuis tensor}
$$N_{\tilde J}(X_1,X_2) =
[\tilde J X_1,\tilde JX_2] -
\tilde J [\tilde J X_1,X_2] - \tilde J [X_1,\tilde
J X_2]  + \tilde J^2 [X_1,X_2],\leqno{(3.2)}$$
where $X_1, X_2 \in \Gamma~  TM$.
Then the almost complex structure $J$ is
projectable iff
$$
N_{\tilde J} \bigr|_{ \Gamma L \times \Gamma TM} \equiv 0.
$$

This is easily checked, by using the fact that
$N_{\tilde J}$ is a tensor. Indeed, consider
the vector fields $Y \in
\Gamma L, X \in \Gamma TM$, extensions of $Y_p \in T_p \cal F$
and $X_p \in T_p M$ $(p\in M)$;
generality is not affected if we assume $X$ projectable,
which, hereafter, we will denote by $X\in\Gamma_{pr}TM$, and which means
that $\forall Y\in\Gamma L$, $[Y,X]\in\Gamma L$.
Then,
$$N_{\tilde J} (Y_{p},X_{p}) =N_{\tilde J} (Y,X)/_{p} =
- \{ \tilde J \bigl( [Y, \tilde J X] -
\tilde J [Y,X] \bigr)\}/_{p} \leqno{(3.3)}$$
$$= -  J \pi [Y,
\tilde J X]_{p}
= -  J (\stackrel{\circ}{D_Y}  J)\pi X/_{p}.$$
Hence, $N_{\tilde J} (Y,X)=0$
if and only if $(\stackrel{\circ}{D_Y}  J) X=0$ , as stated.

From (3.3), we also notice
that, $\forall Y\in \Gamma L,\forall X\in \Gamma E$,
$N_{\tilde J}(Y,X)$ takes values in $E$.

\par It follows that an almost
hypercomplex structure $H=(I_1,I_2,I_3)$ on $E \approx \nu \cal F$
is projectable iff
one has $N_{\tilde I_\alpha} (Y,X)=0$ for $Y \in \Gamma L, X \in \Gamma
TM$, $\alpha=1,2,3$.

This assertion
can be rephrased by
using a unique tensor $T^{\tilde H}: TM \times TM \rightarrow TM$
defined as follows. Recall that for an
almost hypercomplex
structure $H=(I_1,I_2,I_3)$ on a manifold $M^{4q}$, a {\em structure
tensor} is defined by
$$
T^H=  \frac{1}{6} ~ \sum_{\alpha=1}^3 N_{I_\alpha},\leqno{(3.4)}
$$
the torsion of the {\em Obata connection} on $M$
(\cite{AM}, pp. 239-241).
In our case, the almost hypercomplex structure $H=(I_1,I_2,I_3)$ is only
defined on a
complementary
distribution $E$ of the tangent bundle $L$ of the
$4q$-codimensional foliation $\cal F$.
But, we may take the triple $\tilde H=(\tilde I_1,\tilde I_2,\tilde I_3)$
defined as in (3.1), and
define the {\em structure tensor}
$$T^{\tilde H} = \frac{1}{6}\sum_{\alpha=1}^3 N_{\tilde I_\alpha}.
\leqno{(3.5)}$$

The following formula, where $Y\in\Gamma L,X\in\Gamma TM$,
and $\alpha = 1,2,3$, is a
consequence of (3.5)
$$N_{\tilde I_\alpha}(Y,X)=
\frac{3}{2}\{T^{\tilde H}(Y,X)+\tilde I_\alpha T^{\tilde H}(Y,\tilde I_\alpha
X)\}. \leqno{(3.6)}$$
It follows:
\proclaim 3.2 Proposition.
The almost hypercomplex structure
$H=(I_1,I_2,I_3)$ defined on the transversal
bundle $\nu \cal F$ of the foliation $\cal F$ of
$M^{p+4q}$ is projectable iff
$T^{\tilde H} \bigr |_{T {\cal F} \times TM}$ is zero.
\par
Using the tensor (3.5), we can also show another interesting fact namely,
\proclaim 3.3 Proposition. If the foliation ${\cal F}$ has a projectable,
transversal, almost hypercomplex structure $H$, there exists a
projectable
connection of $\nu{\cal F}$ which preserves the structure $H$.\par
\noindent{\bf Proof.} We recall that a Bott connection $\nabla$ of
$\nu{\cal F}$ is projectable if $\nabla_{X_{1}}X_{2}$ is projectable
$\forall X_{1},X_{2}\in\Gamma_{pr}E$. The stated result will be proven
by writing down analogs of connections defined by Oproiu and Obata.
First, let us define a Bott connection $\nabla^{H}$ i.e.,
$\nabla^{H}_{Y}=\stackrel{\circ}{D}_{Y}$ given by (1.1$'$) ($Y\in \Gamma
L$), by adding the equation \cite{{Op2},{AM}}
$$\nabla^{H}_{X_{1}}X_{2}=
{1 \over{12}} \pi
\sum_{(\alpha,\beta,\gamma)} \bigg (\tilde
I_\alpha [I_\beta X_1,I_\gamma X_2] +
\tilde I_\alpha [I_\beta X_2,I_\gamma X_1] \bigg )  \leqno{(3.7)}$$
$$+{1 \over 6} ~\pi\sum_\alpha \bigg (\tilde I_\alpha
[I_\alpha X_1,X_2] + \tilde I_\alpha [I_\alpha X_2,X_1]
\bigg ) + {1 \over 2} \pi [X_1,X_2],$$
where $X_1,X_{2} \in\Gamma E$, $\sum_{(\alpha,\beta,\gamma)}$
denotes the sum over the cyclic
permutations of $(1,2,3)$, and $\pi:TM
\rightarrow E$ is the natural projection. $\nabla^{H}$
is a projectable connection of $\nu{\cal F}$.
It does not preserve $H$ but, if
we correct (3.7) by defining \cite{{Ob},{AM}}
$$D^{H}_{X_{1}}X_{2}=\nabla^{H}_{X_{1}}X_{2}
+\frac{1}{2}\pi T^{\tilde H}(X_{1},X_{2}),
\hspace{5mm}X_{1},X_{2}\in\Gamma E,\leqno{(3.8)}$$
we get a connection as required by the proposition.
The projectability of the additional term of (3.8) follows from (3.2)
since, if $J$ of (3.2) is projectable then
$\forall X_{1},X_{2}\in\Gamma_{pr}E$, $N_{\tilde J}(X_{1},X_{2})$ has a
projectable transversal part.
Q.e.d.

The connection $D^{H}$ of (3.8) will be called the {\em Bott-Obata
connection}, and for its torsion we get
$$T_{D^{H}}(X_{1},X_{2}):=D^{H}_{X_{1}} X_{2}-D^{H}_{X_{2}} X_{1}
-\pi[X_{1},X_{2}]=\pi T^{\tilde H}(X_{1},X_{2}),$$
$\forall X_{1},X_{2}\in\Gamma E$.

Proposition 3.3 shows that, generally, there are obstructions to the
existence of a projectable, transversal, hypercomplex structure of a
foliation ${\cal F}$. One such obstruction is, of course, the Atiyah class
of ${\cal F}$, since the Atiyah class is the obstruction to the existence
of a projectable, transversal connection \cite{Mol}.
Sometimes, it is also
possible to detect secondary characteristic classes.

Let $\cal F$ be a foliation of codimension $4q$ on $M^{p+4q}$, which
has a projectable, transversal almost
complex structure $I_1$. Then there exist
Bott connections $\nabla$ which preserve $I_1$. Indeed,
for any Bott connection $\nabla_{Y}
I_{1}=0$ for all $Y\in\Gamma L$, and the existence of $\nabla$ with
$\nabla_{X}I_{1}=0$ for all $X\in\Gamma E$ follows in the same way as the
existence of, say, an almost complex connection on an almost complex
manifold. Moreover, if we also choose a Riemannian metric $g$ on $M$
such that $g/_{E}$ is $I_{1}$-Hermitian, we can get $\nabla$ as above
which also satisfies $\nabla_{X}(g/_{E})=0$, $\forall X\in\Gamma E$.

Accordingly, as in the classical Bott vanishing theorem \cite{Bt}, we have:
$$Chern_{2k} (E, I_1) =0 \hspace{5mm}{\rm if}\; k>4q,$$
where $Chern_{2k}$ denotes elements of cohomological
degree $2k$ in the ring generated by the real Chern classes. More exactly
the representative differential forms of these classes in terms of the
curvature forms of $\nabla$ vanish.

Now, assume that $I_1$ can be completed by $I_2$, $I_3$
to a (not necessarily projectable) transversal almost
hypercomplex structure, with an almost hyperhermitian metric $g$.
Then, the odd dimensional Chern classes $c_{2h+1}(E,I_{1})$ vanish,
since their representative differential forms
in terms of the curvature of an almost hyperhermitian (not necessarily
Bott)
connection $D$ are ${\cal C}_{2h+1}(D)=0$ (cf. \cite{KN}, vol. II, p. 304).

Thus, if $2h+1 >4q$, and if the connections $\nabla,D$ are as above, we
have
$${\cal C}_{2h+1} (\nabla) - {\cal C}_{2h+1} (D) =
d (\Delta_{(h)} (\nabla,D))=0,$$
where $\Delta_{(h)}$ are the {\em Bott comparison forms} \cite{Bt},
and we get cohomology classes
$$[\Delta_{(h)} (\nabla,D)] \in H^{4h+1}(M, {\bf R}) \hspace{5mm}(h \geq
2q)$$ which are well defined and independent of the choice of the
connections $\nabla,D$.
These precisely are the secondary classes that we mentioned. They
are obstructions to the existence of a projectable almost
hyperhermitian transversal structure with the given
almost complex component $I_1$ since, if such a structure exists, we may
use equal connections $D=\nabla$, in which case
$\Delta_{(h)} (\nabla,D)=0$.\vspace{1mm}

Next, assume that $Q \subset End~  E$ is
a transversal almost quaternionic structure of a foliation ${\cal F}$.
In order to get a tensorial
criterion for the projectability of $Q$,
we look at the extension $\tilde Q \subset
End (TM)$ of $Q$ defined by extending each
$S\in Q$ to $\tilde S\in End\,TM$ by
$\tilde S/_{L}=0$.
Recall that the structure
tensor $T^Q$ of an almost quaternionic structure $Q$
of a manifold $M^{4q}$ is defined by
$$
T^Q (X_1,X_2) = T^H (X_1,X_2) +
\sum_{\alpha=1}^3 [(\tau_\alpha X_1~  )I_\alpha X_2 -
(\tau_\alpha X_2 ~ )
I_\alpha X_1],
\leqno{(3.9)}$$
where $X_{1},X_{2}\in\Gamma TM$,
$H=(I_{1},I_{2},I_{3})$ is any local basis of $Q$, and
$$\tau_\alpha X = \frac{1}{4q-2} ~tr~[(I_\alpha ~T^H) (X,-)],\hspace{5mm}
X\in\Gamma TM$$
(cf. \cite{AM}, p. 244).
Both $T^{H}$ and $T^{Q}$ are invariant by a change of the local basis $H$
since such a change is via an $SO(3)$-matrix and the sums on $\alpha$ which
enter in the expressions of $T^{H},T^{Q}$ behave like scalar products in ${\bf
R}^{3}$.
\par In our situation, $Q$ is
defined only on the complementary distribution $E$ of $L =
T \cal F$, and a suitable extension $T^{\tilde Q}$ of $T^Q$
(i.e., $T^{\tilde Q}(X_{1},X_{2})$ is given by (3.9) with $T^{H}$
replaced by $\pi T^{\tilde H}$, if $X_{1},X_{2}\in\Gamma E)$
will result from
\proclaim 3.4 Proposition. The transversal
almost quaternionic structure $Q$
of the foliation ${\cal F}$ on $M^{p+4q}$ is projectable iff there exists a
local basis $H$ of $Q$ such that
$$T^{\tilde H} (Y,X) = \sum_{\alpha=1}^3 \kappa_\alpha (Y  )\tilde I_\alpha X
\hspace{5mm}(Y\in\Gamma L,\,X\in\Gamma TM)\leqno{(3.10)}$$
for some leafwise $1$-forms $\kappa_{\alpha}$ on $M$.
\par
\noindent{\bf Proof.}
If (3.10) holds, then, $\forall X\in\Gamma E,\,\forall Y\in\Gamma L$,
(3.3) and (3.6) yield
$$
-  I_\lambda (\stackrel{\circ}{D_Y} I_\lambda) X =
N_{\tilde I_\lambda}(Y,X) = \frac{3}{2} \{T^{\tilde
H}(Y,X) +\tilde I_\lambda T^{\tilde H}(Y,\tilde I_\lambda X)\} =
$$
$$=\frac{3}{2}(\sum_{\alpha=1}^3 \kappa_\alpha( Y)
~ I_\alpha X + I_\lambda \sum_{\alpha=1}^3 \kappa_\alpha( Y) ~
I_\alpha (I_\lambda X))
$$
$(\alpha,\lambda=1,2,3)$, whence
$$
(\stackrel{\circ}{D_Y} I_\lambda) X = \frac{3}{2}
\sum_{\alpha=1}^3 \kappa_\alpha (Y) ( I_\lambda ~
I_\alpha X -I_\alpha~ I_\lambda X).\leqno{(3.11)}
$$
Since the structures $I_{\alpha}$ satisfy the quaternionic identities,
(3.11) shows that $Q$ is a projectable structure  (compare with
(1.2)).

Conversely,
if $Q$ is projectable then,
$\forall Y \in \Gamma L, ~
\forall X \in \Gamma_{pr}E$, (1.2) and (3.3) imply
$$
N_{\tilde I_1}(Y,X) =  -\alpha(Y) I_3 X + \beta (Y) I_2
X ,   $$ and similarly:
$$
N_{\tilde I_2}(Y,X) =  - \alpha(Y) I_3 X - \gamma (Y) I_1 X ,
$$
$$
N_{\tilde I_3}(Y,X) =  \beta(Y) I_2 X - \gamma (Y) I_1 X .
$$
Accordingly, (3.5) yields
$$
T^{\tilde H} (Y,X) = - \frac{1}{3} [\gamma (Y) I_1 -
\beta (Y) I_2 + \alpha (Y) I_3 ] X,
$$
which  is (3.10) for $X\in\Gamma E$. For $X\in \Gamma L$, (3.10) is just
$0=0$.
Q.e.d.

Moreover, by
taking into account that $tr~I_\alpha=0$, we get
$$
\alpha(Y) =  \frac{3}{4q}  ~  tr~\{I_3 T^{\tilde H} (Y,-)\},
\hspace{5mm}
\beta (Y) =  -\frac{3}{4q}  ~  tr~\{I_2 T^{\tilde H} (Y,-)\} ,$$
$$\gamma (Y) =  \frac{3}{4q}  ~  tr~\{I_1 T^{\tilde H} (Y,-)\}, $$
where the missing argument is in $\Gamma E$.

Therefore, the coefficients of (3.10) must be $\alpha,\beta,\gamma$,
and we have to define the extension of
$T^{Q}$
by asking $T^{\tilde
Q}(Y_1,Y_2) =0$ for $Y_1,Y_2 \in \Gamma L$, and
$$T^{\tilde Q} (Y,X) =
T^{\tilde H} (Y,X) + \sum_{\alpha=1}^3 \rho_\alpha (Y  )I_\alpha X,
\leqno{(3.12)}$$
for $Y\in \Gamma L, X\in\Gamma E$,
where
$$\rho_\alpha(Y) = (1/4q) tr ~[(I_\alpha ~T^{\tilde H}) (Y,-)].
\leqno{(3.13)}$$
This $T^{\tilde Q}$ is
independent of the choice of the local basis $H$ for the
same reason $T^{Q}$ was.

Accordingly, we see that Proposition 3.4 is equivalent to
\proclaim 3.5 Proposition. The almost quaternionic structure $Q$, transversal
to the foliation $\cal F$ of $M^{p+4q}$, is projectable iff
$T^{\tilde Q} \bigr |_{T {\cal F} \times TM}$ vanishes.\par
Formula (3.11)
gives a geometric meaning to the 1-forms $\rho_\alpha $ of (3.13)
in the case of a projectable structure $Q$.
Namely, they are local connection forms of $\stackrel{\circ}{D_{}}$
restricted to $Q$.
In particular, if the triple $(I_1,I_2,I_3)$ consists of projectable
structures, one has $\rho_\alpha =0$.\vspace{1mm}

It is also interesting to notice that $T^{\tilde Q}(X_{1},X_{2})$
$(X_{1},X_{2}\in\Gamma E)$ can be related with the torsion of some well
chosen Bott connections. First,
all the $Q$-preserving Bott connections on $E$ are given by
$$ \nabla^Q_{X_1} X_2 =\left\{
\begin{array}{ll}
\pi [X_1,X_2],&
\;\;\;\;{\rm for}~  X_1 \in \Gamma L, \vspace{1mm}\\
\nabla^{Op}_{X_1} X_2&
\;\;\;\;{\rm for}~  X_1 \in \Gamma E,
\end{array} \right.\leqno{(3.14)}$$
where $\nabla^{Op}$ denotes
the connection defined by Oproiu's formula (\cite{O}, p. 295)
$$\nabla^{Op}_{X_1} X_2 =\nabla_{X_1} X_2 + \sum_{\alpha=1}^{3} \{\frac{1}{4}
(\nabla_{X_1} I_\alpha) I_\alpha +
\frac{1}{2} \eta_\alpha (X_{1}) I_\alpha \} X_2
\leqno{(3.15)}$$
$$+\frac{1}{4} \{A_{X_1} X_2 - \sum_\alpha I_\alpha A_{X_1}
(I_\alpha X_2 )\}\hspace{5mm}(X_{1},X_{2}\in\Gamma E).$$
In (3.15)
$\nabla$ is an arbitrary Bott connection on $E$, $H=(I_1,I_2,I_3)$ is a
local compatible almost hypercomplex structure,
$A_{X_1}$ is an arbitrary endomorphism of $E$, and $\eta_{\alpha}$
$(\alpha=1,2,3)$ are arbitrary $1$-forms on $M$.

Now, let us fix a connection $\stackrel{1}{\nabla}$ among those given
by (3.14), (3.15).
Following \cite{AM}, p. 244, $\stackrel{1}{\nabla}$
has an {\em associated
Bott-Oproiu connection}
$$^{Op}\stackrel{1}{\nabla}_{X}=\stackrel{1}{\nabla}_{X}
+\sum_{\alpha=1}^{3}(\varphi_{\alpha}+\frac{1}{3}\varphi\circ
I_{\alpha})(X)I_{\alpha}-\frac{1}{4}(A_{X}-\sum_{\alpha=1}^{3}I_{\alpha}A_{X}
I_{\alpha}),\leqno{(3.16)}$$
where
$$\varphi_{\alpha}(X)=\frac{1}{4q-2}tr(I_{\alpha}T_{X}),\;
\varphi=\sum_{\alpha=1}^{3}\varphi_{\alpha}\circ I_{\alpha},\;
A_{X}=T_{X}+\frac{1}{3}\sum_{\alpha=1}^{3}T_{I_{\alpha}X}\circ I_{\alpha},\;
X\in\Gamma M,$$
$T$ being the torsion of $\stackrel{1}{\nabla}$, and $T_{X}$ the
endomorphism of $E$ obtained by fixing the first argument of the torsion as
$X$.

Then, the same computations as in \cite{AM} show that the tensor $T^{\tilde
Q}$ and the torsion of the Bott-Oproiu connections are related by the
formula
$$T_{{^Op}\stackrel{1}{\nabla}}(X_{1},X_{2})=\pi T^{\tilde Q}(X_{1},X_{2}),
\hspace{5mm}X_{1},X_{2}\in\Gamma E,$$
which is the result we wanted to mention.
\par Finally, let us also note that, as a consequence of (3.9), if $Q$ is a
projectable structure, $\pi T^{\tilde Q}$ is a projectable tensor field.
\section{Integrability}
Consider an almost hypercomplex structure
$H=(I_1,I_2,I_3)$, respectively an almost quaternionic structure $Q$
transversal to a foliation
${\cal F}$ of codimension $4q$ on a manifold $M^{p+4q}$.
By Definition 1.2, the integrability of $H$ and $Q$ includes projectability.
It is natural to ask whether integrability can be
recognized by means of the structure tensors $T^{\tilde
H}$ and $T^{\tilde Q}$, defined by formulas (3.5) and (3.9).
\proclaim 4.1 Proposition. $H$, respectively $Q$, is integrable iff
its structure tensor $T^{\tilde H}$, respectively $T^{\tilde Q}$,
takes values in
the tangent bundle $L=T{\cal F}$.
\par
\noindent{\bf Proof.}
If $H$ (respectively $Q$)is projectable, as seen in Section 3,
$\pi\circ T^{\tilde H}$ (respectively $\pi\circ T^{\tilde Q}$)
projects to the local slice spaces, and, clearly,
the projection is the torsion tensor of the corresponding almost
hypercomplex (quaternionic) structures of these slice spaces.
Accordingly, the statement follows by Definition 1.2 and by the fact
that $T^H=0$ (respectively $T^Q=0$)
is the integrability condition for $H$ (respectively $Q$) on manifolds. Q.e.d.

Now, we will discuss another aspect concerning transversal
quaternionic (i.e., integrable, almost quaternionic) structures $Q$ of a
foliation.
The integrability of $Q$ is equivalent to the existence of an open covering
$M=\cup_{a\in{\cal A}}U_{a}$ (${\cal A}$ is an arbitrary set) such that
one has local, torsionless, projectable connections
$D^a$ of $E/_{U_{a}}$
which preserve $Q/_{U_{a}}$. These connections can be glued
together by means of a partition of unity. The resulting global
connection is then a torsionless, $Q$-preserving, Bott connection
but, generally, it is not projectable.
As a matter of fact,
we already have explicit
expressions of such connections namely, the Bott-Oproiu connection of any
$Q$-preserving Bott connection has a vanishing torsion because of
Proposition 4.1. We write this result as
\proclaim 4.2 Proposition. If the foliation $\cal F$ admits a transversal
quaternionic structure $Q$, then ${\cal F}$ admits
a $Q$-preserving, torsionless, Bott connection $D$ on its transversal bundle.
\par
On the other hand, from the system of
local connections $D^a$ above, we can build a \v{C}ech $1$-cocycle,
as follows. The differences
$$ \tau^{ab}(X,Y) = D^a_{X}Y - D^b_{X}Y\hspace{5mm}(X,Y\in\Gamma
E,\,a,b\in{\cal A})\leqno{(4.1)}$$
are projectable cross sections of the foliated vector bundle $Hom(E\odot E,E)$,
where $\odot$ denotes the symmetrized tensor product. Symmetry comes from
the fact that the connections $D^{a}$ have no torsion. Furthermore,
if we denote
$$\tau_{X}^{ab}=D^a_{X} - D^b_{X},\hspace{5mm}X\in\Gamma TM,\leqno{(4.2)}$$
we obtain $(End\,E)$-valued $1$-forms $\tau^{ab}
\in\Lambda^{1}(U_{a}\cap U_{b},End\,E)$, and $\tau^{ab}_{X}(Q)
\subseteq Q$.

Let us denote by $End_{Q}E\subseteq End\,E$ the subbundle of $Q$-preserving
endomorphisms, and
notice the injections of vector bundles
$$i:Hom(E\odot E,E)\rightarrow Hom(TM\otimes TM,E),$$
$$j:\Lambda^{1}(M,End\,E)\rightarrow Hom(TM\otimes TM,E),$$
where $i$ extends a tensor defined on arguments in $E$ to one with
arguments in $TM$ by giving it the value $0$ if an argument is in $L$, and
$$j(\lambda)(X_{1},X_{2}):=\lambda(X_{1})\pi
X_{2},\hspace{5mm}X_{1},X_{2}\in\Gamma TM.$$
The integrability of $Q$ implies that the
$(End\,E)$-valued $1$-forms $\tau$ defined by (4.2) are
projectable cross sections of
the vector bundle
$j^{-1}(i(Hom (E\odot E,E))$ over
$U_{a}\cap U_{b}.$
Thus, the forms $\tau^{ab}$ may be seen as a $1$-cocycle with values in the
sheaf ${\cal S}$ of germs of projectable cross sections of the vector
bundle $j^{-1}(i(Hom (E\odot E,E))$ on $M$.
Of course, ${\cal S}$ is a subsheaf
of germs of projectable $(End\,E)$-valued $1$-forms on $M$.

Correspondingly, we have a cohomology class
$[\tau]_{{\cal S}}\in H^{1}(M,{\cal S})$
associated with the structure $Q$, which we
call the {\em integrability class} of $Q$.

The integrability class can be handled as follows.
The splitting $TM=E\oplus L$ yields a natural
bigrading, called ${\cal F}$-type, of the spaces of
vector fields and differential forms (our convention is to
write the $E$-degree first), and a decomposition of the exterior
differential
$$d=d'_{(1,0)}+d''_{(0,1)}+\partial_{(2,-1)},\leqno{(4.3)}$$
where the indices denote the type of the operators, and $d''$ is
differentiation along the leaves of ${\cal F}$ \cite{{V1},{V4}}.
Following the de Rham type Theorem 4 of \cite{V4},
p.217, $[\tau]_{{\cal S}}$ is the $d''$-cohomology class of a $(1,0)$-form
with values in $j^{-1}(i(Hom (E\odot E,E))$. Namely, put
$$\tau^{ab}=\tau^{a}-\tau^{b},\leqno{(4.4)}$$
where $$\tau^{a}\in\
\Lambda^{1,0}(U^{a},End\,E)\cap\Gamma
(j^{-1}(i(Hom (E\odot E,E))/_{U^{a}}).$$
Then, since $\tau^{ab}$ are projectable forms, the local forms $d''\tau^{a}$
glue up to a global $d''$-closed form ${\cal T}$, and this is the
required representative form of $[\tau]_{{\cal S}}$.

\proclaim 4.3 Proposition. Let $Q$ be a transversal
quaternionic structure of the foliation ${\cal F}$. Then, a
torsionless, projectable,
transversal connection of ${\cal F}$ which preserves $Q$ exists iff
$[\tau]_{{\cal S}}=0$ i.e., iff ${\cal T}$ is $d''$-exact. \par
\noindent{\bf Proof.} $[\tau]=0$ iff one can get relations (4.4) where the
local forms $\tau^{a},\tau^{b}$ are projectable. If this happens, the
operator $$D=D^{a}-\tau^{a}\hspace{5mm}(a\in{\cal A})\leqno{(4.5)}$$
yield a global projectable connection on $E$
which preserves $Q$. Conversely, if $D$ exists, $\tau^{a}=D^{a}-D$ are
projectable, and satisfy (4.4). Q.e.d.

Notice that the (possibly non projectable) connection $D$ of (4.5) exists
for any integrable structure $Q$.
From (4.5) and the projectability of $D^{a}$
it follows that ${\cal T}$ is the $(1,1)$-part of the curvature form of $D$
hence, ${\cal T}$ also represents the Atiyah class of ${\cal F}$ \cite{Mol}.
This proves
\proclaim 4.4 Proposition. The Atiyah class of a transversally quaternionic
foliation belongs to $\iota^{*}(H^{1}(M,{\cal S}))$, where $\iota$ is the
inclusion of ${\cal S}$ into the sheaf of germs of projectable $1$-forms
with values in $End\,E$.\par
In view of the above results, the following terminology is natural.
A projectable, almost quaternionic transversal structure of a foliation
will be called {\em semi-integrable} if it is preserved by a
global, torsionless Bott connection, and it will be called {\em strongly
integrable} if it is preserved by a
global, torsionless, projectable, Bott connection.
\section{The transversal twistor space of $({\cal F},Q)$}
Let $\cal F$ be a foliation of codimension
$4q$ on the manifold $M^{p+4q}$, endowed with a projectable almost
quaternionic structure $Q$ with local bases $(I_{1},I_{2},I_{3})$ on
the transversal bundle $\nu{\cal F} = TM/T\cal F$, and let
$TM = E \oplus L$ ($L=T{\cal F}$) be a chosen splitting, allowing us to
transfer structures between $\nu{\cal F}$ and $E$.
\par Similarly to the case of quaternionic manifolds,
we define the {\em transversal twistor space} of $\cal
F$ by:
$$Z{\cal F} = \{ J \in Q,\quad
J=\alpha_1 I_1 + \alpha_2 I_2 + \alpha_3 I_3,
\quad \alpha_1^2 + \alpha_2^2 +
\alpha_3^2 =1 \},\leqno{(5.1)}$$
i.e.,  $Z{\cal F}$ is the sphere bundle associated with
the Euclidean vector bundle $Q$, where
the metric of $Q$ is that which makes the compatible almost
hypercomplex structures $(I_1,I_2,I_3)$ orthonormal bases.
\par The quaternionic structure $Q$ reduces the structure group of
$E$ to $Gl(q, {\bf H}) \cdot Sp(1)$, and there
exists a corresponding principal bundle $\pi:{\cal B} (E,Q)
\rightarrow M$
of {\em quaternionic frames (bases)}. A frame
$b \in {\cal B} (E,Q)$ may be identified
with an isomorphism $B: ({\bf R}^{4q}, I^o_1,
I^o_2,I^o_3) \rightarrow E$, where the left hand side is equivalent
to the left quaternionic space ${\bf H}^q$, such
that:
$$B^{-1} \circ H\circ B = H^{0}\cdot A,\hspace{5mm}
H=
\left(
\begin{matrix}
I_1\:\\
I_2\:\\
I_3
\end{matrix}
\right),\hspace{2mm}
H^{0}=
\left(
\begin{matrix}
I^o_1\:\\
I^o_2\:\\
I^o_3
\end{matrix}
\right).\leqno{(5.2)}$$
In (5.2), $H$
is an arbitrary almost hypercomplex local basis of $Q$
seen as a line matrix, $H^{0}$ is the canonical basis of ${\bf H}^{q}$
seen as a line matrix,
the composition $\circ$ is for each element of the
line, dot is matrix multiplication, and $A \in SO(3)$.
Accordingly, we may see a quaternionic frame as
$$b= (b_i, b_{i'} = I_1 b_i, b_{i*} = I_2 b_i, b_{i'*} = I_3 b_i)_{i=1}^q
\leqno{(5.3)}$$
where $(b_i)$ is the image by $B$ of the canonical basis of ${\bf H}^q$
over $\bf H$.
\par From formula (5.2) we see that the structure group of the
principal bundle ${\cal
B} (E,Q)$ appears as
$$Gl(q, {\bf H}) \cdot Sp(1) \approx \{\phi \in Aut
({\bf R}^{4q})/ \; \phi^{-1} \circ H^{0}\circ \phi = H^{0}\cdot A,
\; A \in SO(3) \},\leqno{(5.4)}$$
and the corresponding Lie algebra $gl(q, {\bf H}) \oplus sp(1)$ is isomorphic to 
$$\{\chi \in End ({\bf R}^{4q})
/\;H^{0}\circ \chi - \chi \circ H^{0}=
H^{0}\cdot\alpha,
\; \alpha \in so(3) \},\leqno{(5.5)}$$
(cf. \cite{Pi}, p. 595).\vspace{2mm}
\par
For further use, we notice that the dual coframe of $b$ is of the form
$$\beta= (\beta^i, \beta^{i'} = - \beta^i \circ I_1 ,
\beta^{i*} = - \beta^i \circ I_2, \beta^{i'*} = - \beta^i \circ I_3
)_{i=1}^q\leqno{(5.6)}$$
where $\beta^i(b_j) = \delta^i_j$.

Then, $b$ provides the complex frame
$(b_i,b_{i*})$ of $(E,I_{1})$ with the dual coframe $(\beta^i,
\beta^{i*})$. As a complex vector bundle, $(E,I_{1})$ is
isomorphic to the holomorphic part of $E \otimes {\bf C}$, and it is
well known that the corresponding basis of this holomorphic part is
$$c_i = I_1^+ b_i, \quad c_{i*} = I_1^+ b_{i*}\leqno{(5.7)}$$
where:
$$I_1^{+} = \frac{1}{2} (Id - {\sqrt -1} I_1).
\leqno{(5.8)}$$
The dual complex cobasis is:
$$\gamma^i = \beta^i + {\sqrt -1} \beta^{i'} ,
\quad \gamma^{i*} = \beta^{i*} + {\sqrt -1} \beta^{i'*}.
\leqno{(5.9)}$$
\proclaim 5.1 Proposition. ${\cal
B} (E,Q)$ is a foliated principal bundle over $(M,{\cal F})$.
\par
\noindent{\bf Proof.} A foliated structure on a principal bundle is a
maximal local
trivialization atlas with projectable transition functions
e.g., \cite{{Mol},{V3}}.
Consider real local bases of $E$ which have projectable transition functions.
Then, there exists local bases of $E$ over $\bf H$ which consist of
some of the vectors of the given bases, and their images by
the operators
$(I_1,I_2,I_3)$ which span $Q$.
Clearly, if we choose a projectable triple
$(I_1,I_2,I_3)$ (which is possible because of the projectability of $Q$)
the corresponding $\bf H$-bases will also have projectable transition functions.
Q.e.d.
\par Now, from formulas (5.3) and (5.6) it follows that
$$I_1 = \sum_{i=1}^q (b_{i'} \otimes \beta^i - b_i
\otimes \beta^{i'} - b_{i*} \otimes \beta^{i'*} + b_{i'*} \otimes
\beta^{i*}),$$ \vspace{-4mm}
$$ I_2 = \sum_{i=1}^q (b_{i*} \otimes \beta^i - b_{i'*}
\otimes \beta^{i'} - b_{i} \otimes \beta^{i*} + b_{i'} \otimes
\beta^{i'*}),\leqno{(5.10)}$$
$$I_3 = \sum_{i=1}^q (b_{i'*} \otimes \beta^i + b_{i*} \otimes \beta^{i'} -
b_{i'} \otimes \beta^{i*} - b_{i} \otimes
\beta^{i'*}),$$
and these formulas define a projection:
$$\pi_Q : {\cal B} (E,Q) \rightarrow {\cal B} (Q),
\leqno{(5.11)}$$
where ${\cal B} (Q)$ is the $SO(3)$ principal bundle of the
positive orthonormal bases of $Q$.
\par Furthermore, we may also consider the projection
$$\pi_Z : {\cal B} (Q) \rightarrow Z{\cal F}
\leqno{(5.12)}$$
defined by $\pi_Z (I_1,I_2,I_3) = I_1$.
\par Clearly, $\pi_Q$ is a principal fibration with
structure group $GL(q,{\bf H})$, $\pi_Z$ is a principal circle
bundle, and
$\pi_{M}:Z{\cal F} \rightarrow M$
is an associated bundle of ${\cal B} (Q) \rightarrow M$
with group $SO(3)$, and fiber $SO(3)/SO(2) = S^2$.

Furthermore, as a consequence of Proposition 5.1 we have
\proclaim 5.2 Corollary.
There exists a lift ${\widetilde {\cal F}}$ of
$\cal F$ to ${\cal B} (E,Q)$, and $\pi_Z\circ\pi_{Q}$ maps
$\tilde{\cal F}$ onto a
foliation ${\widehat {\cal F}}$ of the
twistor space $Z{\cal F}$.
The leaves of $\tilde{\cal F}$ and
${\widehat {\cal F}}$ are covering spaces of the leaves of $\cal
F$.\par
\noindent{\bf Proof.} A slice of $\tilde{\cal F}$ through
$b\in{\cal B}(E,Q)$ appears as the result of the translation of $b$
along a slice of ${\cal F}$ through $\pi(b)\in M$ by the linear holonomy of
${\cal F}$. And, a slice of $\widehat{\cal F}$ is the result of the
projection of the previous slice of $\tilde{\cal F}$ by
$\pi_{Z}\circ\pi_{Q}$ \cite{Mol}. Q.e.d.

In what follows we will derive local tangent cobases of the
manifold $Z{\cal F}$. We begin by looking
at the ($GL(p,{\bf R})$ $\times$ $(GL(q,{\bf H}) \cdot Sp(1))$)-
principal bundle ${\cal B}(M, Q)$,
consisting of all the tangent bases of $M$ which are of the form
$(a,b)$, $a$ being a
frame of $L$, and $b$ a frame of the form (5.3) in $E$.
The mapping $(a,b)
\rightarrow b$ is a $GL(p,{\bf R})$-principal fibration
$$\pi_{\cal B}:{\cal B}(M, Q) \rightarrow {\cal B}(E, Q).$$
\par On ${\cal B}(M, Q)$, there exists the canonical 1-form \cite{KN}
which, in our
case, has the scalar components, say
$$\alpha^u, \beta^i, \beta^{i'}, \beta^{i*}, \beta^{i'*},
\leqno{(5.13)}$$
where $u=1,...,p$; $i=1,...,q$, and the
forms $\beta$ are as in formula (5.6). From the
known condition \cite{KN}:
$$R_g^* \left(\begin{array}{c}
\alpha\\ \beta\end{array}\right) = g^{-1} \circ
\left(\begin{array}{c}\alpha\\ \beta\end{array}\right),
\leqno{(5.14)}$$
where $\alpha,\beta$ are the columns with the entries defined by (5.13),
and $g \in GL(p,{\bf R}) \times (GL(q,{\bf H}) \cdot Sp(1))$,
it easily follows that the
pullbacks of the forms $\beta$ by local cross sections
of $\pi_{\cal B}$ are global 1-forms on
${\cal B}(E, Q)$ (the transversal canonical 1-form, see \cite{Mol}),
while the pullbacks of $\alpha^u$ yield
some local 1-forms. In this paper,
the pulling back sections will not be written explicitly.
Overall, we get $p+4q$ independent horizontal
(i.e., vanishing on the fibers) 1-forms on ${\cal B}(E, Q)$

Formula (5.14) implies that for any $g \in GL(q,{\bf H}) \cdot Sp(1)$,
and for the
corresponding right translation of the principal bundle
${\cal B}(E, Q)$, one has
$$R_g^* \beta = g^{-1} \circ \beta.\leqno{(5.15)}$$
\par In particular, if $g \in GL(q,{\bf H})$, the
${\bf H}$-version of formula (5.15)
yields right translation formulas of $(\beta^i)$, $(\beta^{i'})$,
$(\beta^{i*})$, $(\beta^{i'*})$
separately. Accordingly, if the forms $\beta$ are pulled back
by local cross sections of
$\pi_{Q}$, one gets local $1$-forms on ${\cal B}(Q)$
such that each of the four
sets of forms above has transition relations of its own i.e., the
annihilator of each set is invariant.
These pullbacks, and those of $(\alpha^u)$ yield $p+4q$
independent {\em horizontal} \cite{KN} local 1-forms on ${\cal B}(Q)$.
\par Then, the same forms
will be pulled back to $Z{\cal F}$ by local cross sections of
$\pi_Z$. Since the composition $\pi_Z \circ \pi_Q$
has right translations which only preserve
the complex structure $I_1$, the 1-forms
obtained in the end on $Z{\cal F}$ have right translation equations
which only preserve the annihilator of
the sets $\{\gamma^i,\gamma^{i*}\}$,
$\{\bar\gamma^i,\bar\gamma^{i*}\} $, defined by formula (5.9).

Finally, after we make a choice of $E$, the column of the forms
$\alpha^u$ also has an invariant annihilator.

The continuation of the building of
nice cobases on $Z{\cal F}$ is by fixing a $Q$-preserving
Bott connection $D$ defined by a 1-form $\varpi$ with values in the Lie
algebra (5.5) on ${\cal B}(E,Q)$. Then, $\varpi$ induces an
$so(3)$-valued connection form $\omega$ on ${\cal B}(Q)$ by
means of the relation:
$$H^{0} \circ \varpi - \varpi \circ H^{0} = H^{0}\cdot\omega.
\leqno{(5.16)}$$
Of
course, both $\varpi$ and $\omega$ vanish on the leaves of the
lifted foliations of ${\cal F}$ to ${\cal B}(E,Q)$ and ${\cal B}(Q)$,
respectively.
Since we see $Z{\cal F}$ as a quotient of ${\cal B}(Q)$, it is the
form $\omega$ which will be of interest.

The 1-forms
$\alpha, \beta, \omega$ provide local tangent cobases on
${\cal B}(Q)$, and if we look at the symmetric
decomposition
$$so(3) =so(2) + m, \qquad \omega = \phi + \psi,
\leqno{(5.17)}$$
where
$$\omega = \left(
\begin{array}{rrr}
0&a&b\\
-a&0&c\\
-b&-c&0
\end{array}
\right), \; \phi = \left(
\begin{array}{rrr}
0&0&0\\
0&0&c\\
0&-c&0
\end{array}
\right), \; \psi = \left(
\begin{array}{rrr}
0&a&b\\
-a&0&0\\
-b&0&0
\end{array}
\right),
\leqno{(5.18)}$$
we see that $\psi$ is a horizontal form on
the principal fibration ${\cal B}(Q)\rightarrow Z{\cal F}$. Thus:
\proclaim 5.3 Proposition.
The pullbacks of the local 1-forms $$\alpha^u,
\beta^i, \beta^{i'}, \beta^{i*}, \beta^{i'*}, a, b$$ to $Z{\cal F}$ by
local cross sections of $\pi_Z$ are local tangent cobases of
the manifold $Z{\cal F}$. Except for $\alpha^u$, all these
forms are of the ${\widehat {\cal F}}$-type $(1,0)$ and
the system of equations $\alpha^u = 0$ is invariant, and it defines
a complementary subbundle $\widehat E$ of $\widehat L =
T{\widehat {\cal F}}$ in the tangent bundle $TZ{\cal F}$. Moreover,
the following system of equations also are invariant by the
transition functions of these cobases, and define subbundles of
$TZ{\cal F} \otimes_{{\bf R}} {\bf C}$:
$$\alpha^{u}=0, \quad \gamma^i = 0, \quad \gamma^{i*} = 0, \quad \xi:
= a+\sqrt{-1}~ b = 0,
\leqno{({\cal C}_{1})}$$
$$\alpha^{u}=0, \quad \gamma^i = 0, \quad \gamma^{i*} = 0, \quad \bar\xi= 0.
\leqno{({\cal C}_2)}$$
\par
\noindent{\bf Proof.} The only thing which has not yet been proven
is the invariance of the equation $\xi =0$. For this, we recall
the formula \cite{KN}
$$ R^*_{\gamma} \omega = \gamma^{-1} \omega \gamma\hspace{5mm}
\gamma \in SO(3).
\leqno{(5.19)}$$
In particular, if
$$ \gamma = \left(
\begin{array}{rrr}
1&0&0\\
0&\cos\phi&\sin\phi\\
0&-\sin\phi&\cos\phi
\end{array}
\right) \in SO(2), $$
we get
$$R^*_\gamma (a,b) = (a \cos\phi -b \sin \phi, a \sin\phi + b \cos\phi)$$
hence,
$$R^*_\gamma \xi = \xi (\cos\phi + \sqrt{-1} \sin\phi).
\leqno{(5.20)}$$
Q.e.d.

The last part of Proposition 5.3 means that we have
\proclaim 5.4 Theorem. The normal bundle $\nu{\widehat {\cal F}}\approx
\widehat E$ is equipped with two almost complex structures $J_1$, $J_2$
which have ${\cal C}_1$, ${\cal C}_2$, respectively, as
bundles of antiholomorphic vectors.\par
\section{Projectability conditions on $Z{\cal F}$}
In this section we find the conditions which ensure that
the almost complex structures $J_1$, $J_2$
are ${\widehat {\cal
F}}$-projectable structures. It was proven in \cite{V3}
that the projectability
conditions are $d''A^\sigma = 0$ (mod. $A^\sigma$)
where $A^\sigma=0$ are the equations of ${\cal C}_1$ and ${\cal C}_2$,
except for $\alpha^{u}=0$,
respectively,
and $d''$ is the ${\widehat {\cal F}}$-leafwise
differential as fixed by the complementary subbundle
$\widehat E$ (see (4.3)).

From the definition of the canonical form \cite{{KN},{Mol}},
and if we use projectable local bases $(I_{1},I_{2},I_{3})$ of $Q$, it
follows that, on $Z{\cal F}$, the forms $\beta$ of (5.13) and the
corresponding $\gamma$ of (5.9), are $\widehat{\cal F}$-projectable. Hence,
$d''\gamma^{i}=0$, $d''\gamma^{i{*}}=0$, which agrees with the above
mentioned projectability condition.

As a matter of fact, we can write down explicit formulas for the
differentials
$d\gamma^i, d\gamma^{i*}$, and we do so since the formulas will
also be needed
later on.
The required differentials are
given by the torsion-structure equations of $\varpi$, which may be written
on $M$ and, then, lifted to ${\cal B}(E,Q)$ or $Z{\cal F}$.

Let us start with a local basis $(I_{1},I_{2},I_{3})$ of $Q$ where the
induced connection $\omega$ has the equations
$$D I_1 =aI_2 +bI_3, \quad D I_2 =-aI_1 +cI_3, \quad D I_3 = -b I_1 -c I_2.
\leqno{(6.1)}$$
This basis can be used to define the frames of (5.3) whence, we see that
the local equations of $\varpi$ can be written as
$$\begin{array}{lcl}
D b_i &= &\hspace*{4mm}\stackrel{0}{\varpi}^j_i b_j +
\stackrel{1}{\varpi}^j_i b_{j'} + \stackrel{2}{\varpi}^j_i b_{j*}
+ \stackrel{3}{\varpi}^j_i b_{j'*},\vspace{1mm}\\
D b_{i'} &=& -\stackrel{1}{\varpi}^j_i b_j +
\stackrel{0}{\varpi}^j_i b_{j'} - (\stackrel{3}{\varpi}^j_i
-a\delta^{j}_{i}) b_{j*}
+ (\stackrel{2}{\varpi}^j_i+b\delta_{i}^{j}) b_{j'*},\vspace{1mm}\\
D b_{i*} &= &-\stackrel{2}{\varpi}^j_i b_j +
(\stackrel{3}{\varpi}^j_i-a\delta_{i}^{j}) b_{j'} +
\stackrel{0}{\varpi}^j_i b_{j*}
-(\stackrel{1}{\varpi}^j_i-c\delta_{i}^{j}) b_{j'*},\vspace{1mm}\\
D b_{i'*} &= &-\stackrel{3}{\varpi}^j_i b_j -
(\stackrel{2}{\varpi}^j_i+b\delta_{i}^{j}) b_{j'} +
(\stackrel{1}{\varpi}^j_i-c\delta_{i}^{j}) b_{j*} +
\stackrel{0}{\varpi}^j_i b_{j'*}.
\end{array}\leqno{(6.2)}$$

Corresponding to these connection equations, there are classical torsion
structure equations which provide the differentials
$d\beta^{i},d\beta^{i'},d\beta^{i^{*}},d\beta^{i'*}$ \cite{{KN},{V3}},
and these equations give us the required formulas
$$d \gamma^i = \gamma^h \wedge (\stackrel{0}{\varpi}^i_h +
\sqrt{-1}\stackrel{1}{\varpi}^i_h) -
\gamma^{h*} \wedge (\stackrel{2}{\varpi}^i_h -
\sqrt{-1}\stackrel{3}
{\varpi}^i_h)\leqno{(6.3)}$$
$$+ \frac{\sqrt{-1}}{2}( \xi \wedge {\bar \gamma}^{i*}
+\bar\xi\wedge\gamma^{i*})+
\gamma^i \circ T_D,
$$
$$d \gamma^{i*} = \gamma^h \wedge (\stackrel{2}{\varpi}^i_h
+ \sqrt{-1}\stackrel{3}{\varpi}^i_h) +
\gamma^{h*} \wedge (\stackrel{0}{\varpi}^i_h -
\sqrt{-1}\stackrel{1}{\varpi}^i_h)\leqno{(6.4)}$$
$$+\frac{\sqrt{-1}}{2} \xi \wedge(\gamma^{i}- {\bar \gamma}^{i})
-\sqrt{-1}c\wedge\gamma^{i*}+\gamma^{i*} \circ T_D,$$
where $T_{D}$ is the torsion of the connection $\varpi$.

Since $T_{D}$ vanishes if one of its arguments is in $L$, we again see that
$d\gamma^{i}$, $d\gamma^{i*}$ do not contain terms in $\alpha^{u}$.
This is another way to justify the equalities
$d''\gamma^{i}=0$, $d''\gamma^{i*}=0$ i.e., the fact that the forms
$\gamma^{i},\gamma^{i*}$ are $\widehat{\cal F}$-projectable $1$-forms.

Now, we must also compute $d\xi$.
First, the structure equations of $\omega$ on $Q$ are:
$$d\omega + \omega \wedge \omega = \Omega ,
\leqno{(6.5)}$$
where, say,
$$\Omega=\left(
\begin{array}{rrr}
0&\cal A&\cal B\\
-\cal A&0&\cal C\\
-\cal B&-\cal C&0
\end{array}
\right) $$ is the curvature matrix of $\omega$.
The entries of $\Omega$ are defined by
$$da = b \wedge c + {\cal A}, \;\; db =- a \wedge c + {\cal B},
\;\; dc = a \wedge b + {\cal C},
\leqno{(6.6)}$$
whence,
$$d \xi = -\sqrt{-1} \xi \wedge c + ({\cal A} + \sqrt{-1}{\cal B}).
\leqno{(6.7)}$$

The 2-forms $\cal A$, $\cal B$ are related to the curvature operator
$R_D$ of $D$. We could obtain this relation by
differentiating (5.16), but we prefer to proceed as follows.
If $\Phi \in \Gamma End~E$ is seen as a $0$-form with
values in $End~E$, and if we denote by ${\bf D}$ the covariant exterior
differential associated with the connection $\varpi$ of $E$,
it is easy to get (cf. \cite{KMS}, Section 11.15)
$${\bf D}^2 \Phi (X_1,X_2) = [R_D (X_1,X_2), \Phi]:=
R_D (X_1,X_2) \circ \Phi - \Phi \circ R_D (X_1,X_2).
\leqno{(6.8)}$$
By applying this formula to $I_1,I_2,I_3$ and using (6.1)
we get
$$\begin{array}{rcr}
{\cal A} I_2 + {\cal B} I_3 &= &[R_D, I_1],\vspace{1mm} \\
-{\cal A} I_1 + {\cal C} I_3 &= &[R_D, I_2], \vspace{1mm}\\
-{\cal B} I_1 - {\cal C} I_2 &= &[R_D, I_3]
\end{array}\leqno{(6.9)}$$

In order to solve equations (6.9), we use the canonical Euclidean
metric $<\;,\;>_{Q}$ of the $SO(3)$-vector bundle $Q$, while identifying
the Lie algebra $so(3)$ with the Euclidean space ${\bf R}^{3}$. Then,
$<\;,\;>_{Q}$ corresponds to the scalar product, and composition of
endomorphisms, elements of $Q$, to the vector product of vectors of
${\bf R}^{3}$. The solutions are
$$\begin{array}{lll}
\cal A &= &<I_2,[R_D,I_1]>_Q,\vspace{1mm}\\
\cal B &= &<I_3,[R_D,I_1]>_Q =-<I_2,I_1 \circ [R_D,I_1]>_Q,\vspace{1mm}\\
\cal C &= &<I_3,[R_D, I_2]>_{Q}.
\end{array}\leqno{(6.10)}$$

Accordingly, (6.7) becomes
$$d \xi = -\sqrt{-1} \xi \wedge c + <I_2, [R
_D,I_1] + \sqrt{-1} [R_D, I_1] \circ I_1>_Q.
\leqno{(6.11)}$$

Now, the projectability conditions left are
$$\begin{array}{lll}
 d''\xi =0 & ({\rm mod.} ~\gamma^i, \gamma^{i*}, \xi)&
{\rm for} ~J_1 ,\vspace{1mm}\\
 d''\xi =0 &({\rm mod} ~{\bar \gamma}^i, {\bar \gamma}^{i*}, \xi)&
{\rm for}~ J_2.
\end{array}\leqno{(6.12)}$$

With (6.7), the meaning of the projectability conditions (6.12) is
$$({\cal A} + \sqrt{-1} {\cal B})(c_i, Y) =0,\hspace{2mm}
({\cal A} + \sqrt{-1} {\cal B})(c_{i^*}, Y) =0,\hspace{2mm}
{\rm for}\;J_{1},\leqno{(6.13)}$$
\vspace{-8mm}
$$({\cal A} + \sqrt{-1} {\cal B})({\bar c}_i, Y) =0,\hspace{2mm}
({\cal A} + \sqrt{-1}{\cal B})({\bar c}_{i^*}, Y) =0,\hspace{2mm}
{\rm for}\;J_{2},\leqno{(6.14)}$$ correspondingly,
where $(c_{i},c_{i^{*}})$ were defined in (5.7), and $Y\in \Gamma L$.

Since for any vector $X \in \Gamma E$, $I_1^+ X$ can play
the role of $c_i$ for one frame, and of
$c_{i^*}$ for another frame, and $I_{1}^{-}X:=\overline{I_{1}^{+}X}$
can play the role of
$({\bar c}_i,{\bar c}_{i^*})$, respectively,
the projectability
conditions become
$$
({\cal A} + \sqrt{-1} {\cal B})(X - \sqrt{-1} I_1 X, Y)=0,\leqno{(6.15)}$$
\vspace{-8mm}
$$({\cal A}
+ \sqrt{-1} {\cal B})(X + \sqrt{-1} I_1 X, Y) =0,
\leqno{(6.16)}$$
for $J_{1}$ and $J_{2}$, respectively, and where
$X \in \Gamma E$, $Y \in \Gamma L$.
\par If the real and imaginary parts are separated, this means
$$ \begin{array}{r}
<I_2, [R_D (X,Y), I_1] + [R_D (I_1 X,Y),I_1]\circ I_1> =0, \vspace{2mm}\\
<I_2, [R_D (I_1 X,Y), I_1] - [R_D (X,Y),I_1]\circ I_1> =0,
\end{array}\leqno{(6.17)}$$ for $J_{1}$, and
$$\begin{array}{r}
<I_2, [R_D (X,Y), I_1] - [R_D (I_1 X,Y),I_1]\circ I_1> =0, \vspace{2mm}\\
<I_2, [R_D (I_1 X,Y), I_1] + [R_D (X,Y),I_1]\circ I_1> =0,
\end{array}\leqno{(6.18)}$$ for $J_{2}$.

Now, if the basis $(I_1,I_2,I_3)$ of $Q$ is changed to $(I_1, -I_3,I_2)$,
the same conditions will
hold for $I_3$ instead of $I_2$, which means that we
have to replace the projectability conditions of $J_1$ by
$$
[R_{D}(X,Y),I_{1}] + [R_{D}(I_{1}X,Y),I_1]\circ I_{1} = \mu I_{1},
$$  
\vspace{-1.5cm}
$$\leqno{(6.19)}$$
\vspace{-1cm}
$$[R_D (I_1 X,Y), I_1] - [R_D (X,Y),I_1]\circ I_1 = \nu I_1,
$$
and those of $J_2$ by
$$
[R_D (X,Y), I_1] - [R_D (I_1 X,Y),I_1]\circ I_1 =\mu' I_1,$$
\vspace{-1.5cm}
$$\leqno{(6.20)}$$
\vspace{-1cm}
$$[R_D (I_1 X,Y), I_1] + [R_D (X,Y),I_1]\circ I_1 =\nu' I_1,
$$
where, in fact, $I_1$ is any $S \in Q$, $S^2 = -Id$.
Furthermore, if we take the trace in (6.19), (6.20), we get
$\mu = \nu =\mu ' = \nu ' =0.$
Then, in both (6.19) and (6.20), the second
relation is the first
composed by $I_1$. Therefore, the projectability conditions reduce to
$$
[R_D (X,Y), S] + [R_D (S X,Y),S]\circ S =0, \hspace{5mm}
{\rm for}\;J_{1},$$
\vspace{-1.5cm}
$$\leqno{(6.21)}$$
\vspace{-1cm}
$$[R_D (X,Y), S] - [R_D (S X,Y),S]\circ S =0,\hspace{5mm}
{\rm for}\;J_{2}.
$$

Since these conditions are tensorial, it suffices to write them
for a projectable cross
section $S$ of $Q$, and a projectable vector field $X$.
Using $$R_D(X,Y) = [D_X,D_Y] - D_{[X,Y]}$$ we get
$$\begin{array}{ll}
D_Y (D_X S) - S D_Y (D_{S X} S) = 0, &\hspace*{5mm}{\rm
for}\;J_{1},\vspace{2mm}\\
D_Y (D_X S) + S D_Y (D_{S X} S) = 0, &\hspace*{5mm}{\rm for}\;J_{2}.
\end{array}
\leqno{(6.22)}$$
These formulas give us the final form of the projectability conditions:
\proclaim 6.1 Theorem.
(a) the structure
$J_1$ is projectable iff  $\forall X \in \Gamma_{pr} E$ and for
any projectable
cross section $S$ of $Q$
the endomorphism $D_{X}S -
SD_{SX}S$ is projectable.\\
(b) the structure
$J_2$ is projectable iff $\forall X \in \Gamma_{pr} E$ and for any projectable
cross section $S$ of $Q$
the endomorphism $D_X S +
SD_{SX}S$ is projectable.\\
(c) $J_{1}$ and $J_{2}$
are both projectable iff the connection induced by $D$ in $Q$ is
projectable.
\par
\section{Integrability conditions on $Z{\cal F}$}
Now, let us assume that we are in the case where $J_1$, $J_2$ are both
projectable, and study the
integrability of these structures.
\par In this case, and if we use projectable local bases
$(I_{1},I_{2},I_{3})$ of the projectable, transversal, almost quaternionic
structure $Q$,
$\gamma^i, \gamma^{i^*}$ and $\xi$
are ${\widetilde {\cal F}}$-projectable (see (6.3), (6.4) and Theorem
6.1 (c)), and it remains to ask that, for arguments in $E$, one had
$$d\gamma^i =0,\qquad d\gamma^{i^*} = 0,
\qquad d\xi =0 \qquad ({\rm mod.} ~ \gamma^i, \gamma^{i^*}, \xi)
\leqno{(7.1)}$$
for $J_1$, and
$$d\gamma^i =0,\qquad d\gamma^{i^*} = 0, \qquad d{\bar \xi} =0
\qquad ({\rm mod.} ~  \gamma^i,
\gamma^{i^*},
{\bar \xi})
\leqno{(7.2)}$$
for $J_2$.
\par From (6.4), we see that (7.2) never holds.
Thus, $J_{2}$ is never integrable, and we do not have to worry about it
anymore.

Furthermore, (6.3) and (6.4)
yield a {\em torsion integrability condition} of $J_1$ namely,
$$\gamma^i \circ T_D  =0, \qquad \gamma^{i^*} \circ T_D = 0
\qquad ({\rm mod.} \gamma^i, \gamma^{i^*}).
\leqno{(7.3)}$$
The forms (7.3) are the holomorphic components
of $T_D$, i. e., of $I_1^+ \circ T_D$, and (7.3)
means that $I_1^+ \circ T_D$ must vanish on arguments
of the form $$I_1^- X,\hspace{5mm} I_1^- I_2 X = \frac{1}{2}
(I_2 X + \sqrt{-1} I_3 X).$$

If we assume $q \geq 2$, independent arguments
$I_1^- X_1, I_1^- X_2$ exist, and
the torsion integrability condition reduces to
$$I_1^+(T_D(I_1^- X_1, I_1^- X_2))=0,
\leqno{(7.4)}$$
$\forall I_1 \in Q$, $\forall X_1,X_2 \in \Gamma E$.
The explicit form of
(7.4) is
$$\begin{array}{cc}
T_D (X_1 + \sqrt{-1} I_1 X_1, X_2 + \sqrt{-1} I_1 X_2)
\vspace{2mm}\\- \sqrt{-1}
I_1 T_D (X_1 + \sqrt{-1} I_1 X_1, X_2+
\sqrt{-1} I_1 X_2) = 0,
\end{array} \leqno{(7.5)}$$
where, in fact, $I_1$ is any $S\in Q$, $S^{2}=-Id$.
Then, after we separate the real and imaginary part of (7.5), we get
the integrability conditions
$$T_D (X_1,X_2) - T_D (S X_1, S X_2) + S T_D (S X_1, X_2)+
S T_D (X_1, S X_2) =0,
\leqno{(7.6)}$$
$$T_D (S X_1,X_2) + T_D ( X_1, S X_2) - S T_D ( X_1, X_2)+
S T_D (S X_1, S X_2) =0.
\leqno{(7.7)}$$
Since (7.7) is the result of composing (7.6) by $S$ at the left, we get
\proclaim 7.1 Proposition. If $q\geq2$,
the torsion integrability condition of $J_{1}$ is {\rm (7.6)}
$\forall S \in Q$, $S^{2}=-Id$,
and $\forall X_1,X_2 \in \Gamma E$.
In particular, this condition holds if $T_D =0$. \par
Furthermore, we also have a curvature integrability condition
which follows from (6.7) namely, that on arguments in $E$ one had
$${\cal A} + \sqrt{-1} {\cal B} = 0 \hspace{3mm}
({\rm mod.} \;\gamma^i, \gamma^{i^*}).
\leqno{(7.8)}$$

If $q \geq 2$, all we have to ask is that, for all $X_1,X_2 \in \Gamma E$,
the following relation holds:
$$({\cal A} + \sqrt{-1} {\cal B}) (X_1 + \sqrt{-1} I_1 X_1, X_2 +
\sqrt{-1} I_1 X_2) =0.\leqno{(7.9)}$$
The imaginary part of (7.9) is equivalent to its real part by the
transformation $X_1 \mapsto I_1 X_1$. Therefore,
the only remaining curvature integrability condition is
$${\cal A} (X_1,X_2) - {\cal A} (I_1 X_1, I_1 X_2) -
{\cal B} (I_1 X_1, X_2) - {\cal B} (X_1,I_1 X_2) = 0.
\leqno{(7.10)}$$

Here $\cal A$ and $\cal B$ are given by (6.10),
which transforms (7.10) into
$$\begin{array}{c}
<I_2,[R_D(X_1,X_2),I_1] -[R_D(I_1 X_1,I_1 X_2), I_1]\vspace{2mm}\\
+I_1 \circ [R_D (I_1 X_1,X_2), I_1] +
I_1 \circ [R_D (X_1,I_1 X_2), I_1]>_Q =0.
\end{array}
\leqno{(7.11)}$$

Now, note that condition (7.11) must be imposed for any basis
$(I_1,I_2,I_3)$ of $Q$ hence, if
$(I_1,I_2,I_3) \mapsto (I_1,-I_3,I_2)$, we get the same
relation (7.11) for $I_3$ instead
of $I_2$. It follows that the second factor of the scalar product
(7.11) must be of the form
$\lambda I_1$. Then, by taking the trace as we did in (6.19), (6.20), we
get $\lambda=0$.
Moreover, we may take any $S\in Q$, $S^{2}=-Id$ as $I_1$.
\par Therefore, we have obtained
\proclaim 7.2 Theorem.
The curvature integrability condition of $J_1$ is:
$$\begin{array}{c}
[R_D(X_1,X_2),S] -[R_D(S X_1,S X_2), S]\vspace{2mm} \\
+S \circ [R_D (S X_1,X_2), S] + S \circ [R_D (X_1,S X_2), S]=0.
\end{array}\leqno{(7.12)}$$
$\forall S \in Q$ and $\forall X_1,X_2 \in \Gamma E$.\par
\proclaim 7.3 Remark.
Lemma {\rm 14.74} of
{\rm \cite{B}} tells us that, if $T_D =0$, {\rm (7.12)} holds.
In particular, if $Q$ is strongly integrable (see Section 4), and if
$D$ is a $Q$-preserving, projectable, torsionless connection, $J_{1}$ is
integrable. \par
For $q=1$, since we have only one independent
vector $c_1$, which can be obtained from an arbitrary $X$, the
torsion integrability condition is
$$I_1^+ (T_D (I_1^- X, I_2 X + \sqrt{-1} I_3 X))=0,
\leqno{(7.13)}$$
where the real and imaginary parts are equivalent by $X \mapsto I_1 X$.
Hence, (7.13) reduces to
$$T_D (X, I_2 X) - T_D (I_1 X, I_3 X) + I_1 T_D (X,I_3 X) +
I_1 T_D (I_1 X, I_2 X) = 0,
\leqno{(7.14)}$$
which has to hold for any basis $(I_1,I_2,I_3)$ of $Q$.\vspace{2mm}
\par
Furthermore, for $q=1$, the curvature integrability condition is
$$({\cal A} + \sqrt{-1} {\cal B}) (I_1^- X, I_2 X + \sqrt{-1} I_3 X) = 0.
\leqno{(7.15)}$$
and, if we separate the real and imaginary parts, we get
$$\begin{array}{lcl}
{\cal A} (X,I_2 X) - {\cal A} (I_1 X, I_3 X) -
{\cal B} (I_1 X, I_2 X) - {\cal B} (X,I_3 X) &=&
0,
\vspace{1mm}\\
{\cal A} (I_1 X,I_2 X) + {\cal A} (X, I_3 X) +
{\cal B} (X, I_2 X) - {\cal B} (I_1 X,I_3 X) &=&
0.
\end{array}
\leqno{(7.16)}$$
Now, if $(I_1,I_2,I_3) \mapsto (I_1, -I_3,I_2)$, then  $(\cal A,
\cal B) \mapsto (-\cal B, \cal A)$, and the first relation (7.16)
becomes the
second. Hence the only remaining condition is
$$\begin{array}{c}
<I_2,[R_D(X,I_2 X),I_1] -[R_D(I_1 X,I_3 X), I_1]\vspace{2mm}\\
+I_1 \circ [R_D (I_1 X, I_2 X), I_1] + I_1 \circ [R_D (X,I_3 X), I_1]>_Q =0,
\end{array}
\leqno{(7.17)}$$
for all the local, hypercomplex bases of $Q$.
\par If we write equation (7.17) for $(I_1, -I_3, I_2)$,
replacing the first factor $I_2$
by $I_3 \circ I_1$ and using $<\Phi \circ \psi, \chi>_Q =
<\Phi, \psi \circ \chi>_Q$, the result is
again (7.17), where the first factor is replaced by $I_3$.
Hence, the second factor of the scalar product is
proportional to $I_1$, and using the trace as we already did,
this second factor must be zero. Therefore, the
curvature integrability condition becomes
$$\begin{array}{c}
[R_D(X,I_2 X),I_1] -[R_D(I_1 X,I_3 X), I_1]\vspace{2mm}\\
+I_1 \circ [R_D (I_1 X, I_2 X), I_1] + I_1 \circ [R_D (X,I_3 X), I_1] =0,
\end{array}\leqno{(7.18)}$$
for any orthonormal basis of $Q$ and $\forall X \in \Gamma E$.
\par The case $q=1$ is that of a four-dimensional conformal structure.
Hence, the obtained
conditions must be equivalent with those of classical twistor theory.

Coming back to the torsion integrability condition of $J_{1}$, we will
notice the following interesting fact
\proclaim 7.4 Proposition.
Two $Q$-preserving Bott connections $\varpi,\varpi'$
define the same structure $J_1$ on $Z(\cal F)$ iff
their torsions differ by a term which satisfies the
torsion integrability condition.\par
\noindent{\bf Proof.} From the definition of $J_{1}$, it follows that
$\varpi,\varpi'$
define the same structure $J_{1}$ iff the (horizontal) difference form of
the connections induced in $Q$ satisfies
$$\xi'-\xi=0\;\;({\rm mod.}\gamma^{i},\gamma^{i^{*}})\leqno{(7.19)}$$

By subtracting the corresponding structure equations
(6.3), (6.4) of the two connection forms $\varpi,\varpi'$, we get
$$\gamma^h \wedge A^i_h  - \gamma^{h^*} \wedge B^i_h
+\frac{\sqrt{-1}}{2}[(\bar\xi'-\bar\xi)\wedge\gamma^{i*}+ (\xi' - \xi)
\wedge {\bar \gamma}^{i^*}]\leqno{(7.20)}$$
$$+ \gamma^i \circ (T_{D'} - T_D)=0,$$
$$ - \gamma^h \wedge C^i_h  - \gamma^{h^*} \wedge S^i_h
+\frac{\sqrt{-1}}{2} (\xi' - \xi)
\wedge (\gamma^{i}-{\bar \gamma}^{i})-\sqrt{-1}(c'-c)\wedge\gamma^{i*}
\leqno{(7.21)}$$
$$+ \gamma^{i^*} \circ (T_{D'} - T_D)=0,$$
where $A, B, C, S$ are the entries of the difference forms of the
connections, and $D,D'$ are the corresponding covariant derivatives.
Hence, $\xi' -
\xi$ may be calculated by applying $i(I_1^- X)$, with $X \in\Gamma E$
to (7.20), (7.21), and the
result contains only $\gamma^i, \gamma^{i^*}$ iff
$(T_{D'} - T_D)$ satisfy (7.4), if $q\geq2$, and (7.13), if $q=1$.
Q.e.d.
\proclaim 7.5 Corollary.
Any two connections which
satisfy the torsion integrability condition
define the same structure $J_1$.
\par
It is also possible to find the condition for two connections
$\varpi,\varpi'$ as in
Proposition 7.4 to define the same
pair of structures $(J_{1},J_{2})$. Namely,
\proclaim 7.6 Proposition. Under the hypotheses of Proposition 7.4, the
connections $\varpi,\varpi'$ define the same structures $J_{1},J_{2}$ iff one of the
following equivalent conditions is satisfied:\\
(a) $\varpi$ and $\varpi'$ induce the same connection in the vector bundle $Q$;\\
(b) $\forall S\in Q$ with $S^2=-Id$,
one has equal commutants $[S,\varpi]=[S,\varpi']$;\\
(c) $\forall S\in Q$
with $S^2=-Id$, one has equal traces $tr(S\varpi)=tr(S\varpi')$.\par
\noindent{\bf Proof.} Now, we ask condition (7.19) and also
$$\xi' - \xi = 0\hspace{3mm} ({\rm mod.} \bar \gamma^i, \bar
\gamma^{i^{*}}),\leqno{(7.22)}$$
which expresses the fact that the two connections define the same structure
$J_2$.  Together, (7.19) and (7.22) yield
$\xi' - \xi = 0$, which exactly is condition (a).

Furthermore, let us look at the relation (5.16) between the connection form
$\varpi$ and the connection form $\omega$ of the connection induced in $Q$.
With the notation (5.18), and like for (6.10), we obtain
$$a = <I_2, [I_1, \varpi]>_Q, \qquad b= <I_3, [I_1, \varpi]_Q>.
\leqno{(7.23)}$$

If $t:= \varpi' - \varpi$, and in view of condition (a),
we will have $(J_1,J_2) = (J'_1,J'_2)$ iff
$$<I_2, [I_1,t]>_Q =0, \qquad <I_3, [I_1,t]>_Q = 0.
\leqno{(7.24)}$$
This implies $[I_1,t] = \alpha I_1$, and the trace yields
$\alpha = 0$, which exactly is condition (b).

Finally, we will obtain condition (c) by using the following
straightforward solutions of (5.16):
$$\begin{array}{lcr}
a~~Id &=& \frac{1}{2} \{I_3 \varpi + \varpi I_3 +
I_2 \varpi I_1 - I_1 \varpi I_2 \},\vspace{2mm} \\
b~~Id &=& -\frac{1}{2} \{I_2 \varpi + \varpi I_2 -
I_3 \varpi I_1 + I_1 \varpi I_3 \}, \vspace{2mm}\\
c~~Id &=& \frac{1}{2} \{I_1 \varpi + \varpi I_1 -
I_2 \varpi I_3 + I_3 \varpi I_2 \}.
\end{array}\leqno{(7.25)}$$
By taking the traces in (7.25), we get
$$a = \frac{1}{2q} ~tr~(I_3 \varpi), \quad b
= -\frac{1}{2q} ~tr~(I_2 \varpi), \quad c = \frac{1}{2q}
~tr~(I_1 \varpi).
\leqno{(7.26)}$$
Obviously, this proves that (c) is equivalent to (a). Q.e.d.
\proclaim 7.7 Remark. It is also possible to solve equations (6.9) by
formulas of the type (7.25):
$$\begin{array}{lcl}
{\cal A} \, Id &=& \frac{1}{2}\{[R_D, I_3]+[R_D,I_2] 
\circ I_1 - [R_D,I_1]\circ I_2\},\vspace{2mm}\\
{\cal B} \, Id &= &\frac{1}{2}\{[R_D, I_2]-
[R_D,I_3] \circ I_1 + [R_D,I_1]\circ I_3\},\vspace{2mm}\\
{\cal C} \, Id &=& \frac{1}{2}\{[R_D, I_1]-
[R_D,I_2] \circ I_3 + [R_D,I_3]\circ I_2\},
\end{array}\leqno{(7.27)}$$
which implies that $\cal A, \cal B, \cal C$ are
$(1/4q)$ of the traces
of the right hand sides of (7.27). These formulas can provide another form of
writing the projectability and integrability conditions.\par

\vspace{1cm}
{\small Dipartimento di Matematica}\\
{\small Universit\`a di Roma "La Sapienza", Italy}\\
{\small e-mail: piccinni@mat.uniroma1.it}\\
\medskip
\par \noindent
{\small Department of Mathematics}\\
{\small University of Haifa, Israel}\\
{\small e-mail: vaisman@math.haifa.ac.il}
\end{document}